\documentclass[11pt,psamsfonts]{amsart}
\usepackage{amsmath}
\usepackage{amsthm}
\usepackage{amssymb}
\usepackage{amscd}
\usepackage{amsfonts}
\usepackage{amsbsy}
\usepackage{epsfig,afterpage}
\usepackage[dvips]{psfrag}
\usepackage[all]{xy}
\usepackage{pstcol}
\newcommand{\R}{\ensuremath{\mathbb{R}}}

\newcommand{\N}{\ensuremath{\mathbb{N}}}

\newcommand{\Z}{\ensuremath{\mathbb{Z}}}

\newcommand{\e}{\epsilon}

\newtheorem {theorem} {Theorem} 
\newtheorem {proposition}  {Proposition}
\newtheorem {corollary}  {Corollary}

\newtheorem {definition}  {Definition}
\newtheorem {remark}  {Remark}

\newcommand{\bbox}{\ \hfill\rule[-1mm]{2mm}{3.2mm}}

\definecolor{red}{rgb}{1.,0.,0.}
\definecolor{blue}{rgb}{0.,0.,1.}
\definecolor{pink}{rgb}{1.,0.75,0.8}

\begin{document}

\title[Canard Cycles] {Canard Cycles and Poincar\'{e} Index of Non-Smooth Vector Fields on the Plane}

\author[ C.A. Buzzi, T. de Carvalho, P.R. da Silva]
{Claudio A. Buzzi$^1$, Tiago de Carvalho$^2$ and Paulo R. da
Silva$^3$}

\address{$^1$ $^2$ $^3$ IBILCE--UNESP, CEP 15054--000
S. J. Rio Preto, S\~ao Paulo, Brazil}

\email{$^1$buzzi@ibilce.unesp.br  $^2$ti-car@hotmail.com }

\email{$^3$prs@ibilce.unesp.br}


\date{}
\dedicatory{} \maketitle


\begin{abstract}
This paper is concerned with closed orbits of non-smooth vector
fields on the plane. For a subclass of non-smooth vector fields we
provide necessary and sufficient conditions for the existence of
canard kind solutions. By means of a regularization we prove that
the canard cycles are singular orbits of singular perturbation
problems which are limit periodic sets of a sequence of limit
cycles. Moreover, we generalize the Poincar\'{e} Index for
non-smooth vector fields.

\end{abstract}


\section{Introduction}

 Piecewise-smooth systems are widespread within application
areas such as engineering, economics, medicine, biology and ecology.
 The most common piecewise-smooth
systems involve either a discontinuity in the vector field, or in
the orbit given by the integral solution $x(t)$. In this paper we
consider the former, that is, general systems where the vector field
is independently defined on either side of a smooth codimension one
switching manifold. Three possible regions of the manifold are then
apparent. At a crossing region the component of the vector field
normal to the switching manifold has the same direction on both
sides of the manifold (sometimes called sewing instead of crossing).
At a stable sliding region both normal components of the vector
field point toward the manifold. At an unstable sliding region both
normal components point away from the manifold. Piecewise-smooth
systems with sliding are also known as Filippov systems. Clearly
these three different scenarios lead to vastly different dynamics.
An orbit that meets the switching manifold at a crossing region
passes through it, but is non-differentiable at the crossing point.
An orbit that impacts at a stable sliding region sticks becomes
constrained (sticks) to the manifold. An orbit in an unstable
sliding region slides along the switching manifold, but will depart
it under any infinitesimal perturbation. Consequently, the only
means by which a stable sliding orbit can escaping the switching
manifold is tangentially, at the boundary of the sliding region.
This leads to the observation that, under parameter variation,
orbits in Filippov systems can undergo a large variety of
bifurcations, commonly called sliding bifurcations.\\

In this paper we study piecewise-smooth system  on open regions on
the plane. Let $\mathcal{U} \subseteq \R ^{2}$ be an open set and
$\Sigma \subseteq \mathcal{U}$ given by $\Sigma =f^{-1}(0),$ where
 $f:\mathcal{U} \longrightarrow \R$ is a smooth function having $0\in
\R$ as a regular value (i.e. $\nabla f(p)\neq 0$, for any $p\in
f^{-1}({0}))$. Clearly $\Sigma$ is the separating boundary of the
regions $\Sigma_+=\{q\in \mathcal{U} | f(q) \geq 0\}$ and
$\Sigma_-=\{q \in \mathcal{U} | f(q)\leq 0\}$. We can assume that
$\Sigma$ is represented, locally
around a point $q=(x,y)$, by the function $f(x,y)=y.$\\

Designate by $\chi^r$ the space of $C^r$ vector fields on
$\mathcal{U}$ endowed with the $C^r$-topology with $r\geqslant1$ or
$r=\infty,$ large enough for our purposes. Call
\textbf{$\Omega^r=\Omega^r(\mathcal{U},f)$} the space of vector
fields $X_{0}: \mathcal{U} \setminus\Sigma \longrightarrow \R ^{2}$
such that

\vspace{-.5cm}

\begin{equation}
\label{discontinuity} X_{0}(x,y)=\left\{\begin{array}{l}
X_1(x,y),\quad $for$ \quad (x,y) \in \Sigma_+,\\ X_2(x,y),\quad
$for$ \quad (x,y) \in \Sigma_-,
\end{array}\right.
\vspace{-.2cm}\end{equation} where $X_i=(f_i,g_i)\in \chi^r, i=1,2.$
We write $X_0=(X_1,X_2),$ which we will accept to be multivalued in
the points of $\Sigma.$ The trajectories of $X_0$ are solutions of
 ${\dot q}=X_0(q),$ which
has, in general, discontinuous righthand side. The basic results of
differential equations, in this context, were stated by Filippov in
\cite{Fi}. Related theories can be found in \cite{K, ST, T}. In this
paper we consider finite discontinuities, i.e., the vector fields
$X_1$ and $X_2$ are defined in the set $f^{-1}(0)$. Another kind of
discontinuity of which the vector field tends to infinity when it
approximates to the switching manifold
 can be obtained in the equations with impasse (more details in \cite{SZ}).\\

In what follows we will use the notation $ X_i.f(p)=\left\langle
\nabla f(p), X_i(p)\right\rangle $.

We  distinguish the following regions on the discontinuity set
$\Sigma:$
\begin{itemize}
\item [(i)]$\Sigma_1\subseteq\Sigma$ is the \textit{sewing region} if
$(X_1.f)(X_2.f)>0$ on $\Sigma_1$ .
\item [(ii)]$\Sigma_2\subseteq\Sigma$ is the \textit{escaping region} if
$(X_1.f)>0$ and $(X_2.f)<0$ on $\Sigma_2$.
\item [(iii)]$\Sigma_3\subseteq\Sigma$ is the \textit{sliding region} if
$(X_1.f)<0$ and $(X_2.f)>0$ on $\Sigma_3$.
\end{itemize}

\begin{figure}[!h]
\begin{minipage}[b]{0.485\linewidth}
\begin{center}
\psfrag{A}{$q$} \psfrag{B}{$X_0^\Sigma(q)$} \psfrag{C}{$\Sigma_3$}
\psfrag{D}{$q + X_1(q)$} \psfrag{E}{$q + X_{2}(q)$} \psfrag{F}{$m$}
\psfrag{G}{} \epsfxsize=5.5cm \epsfbox{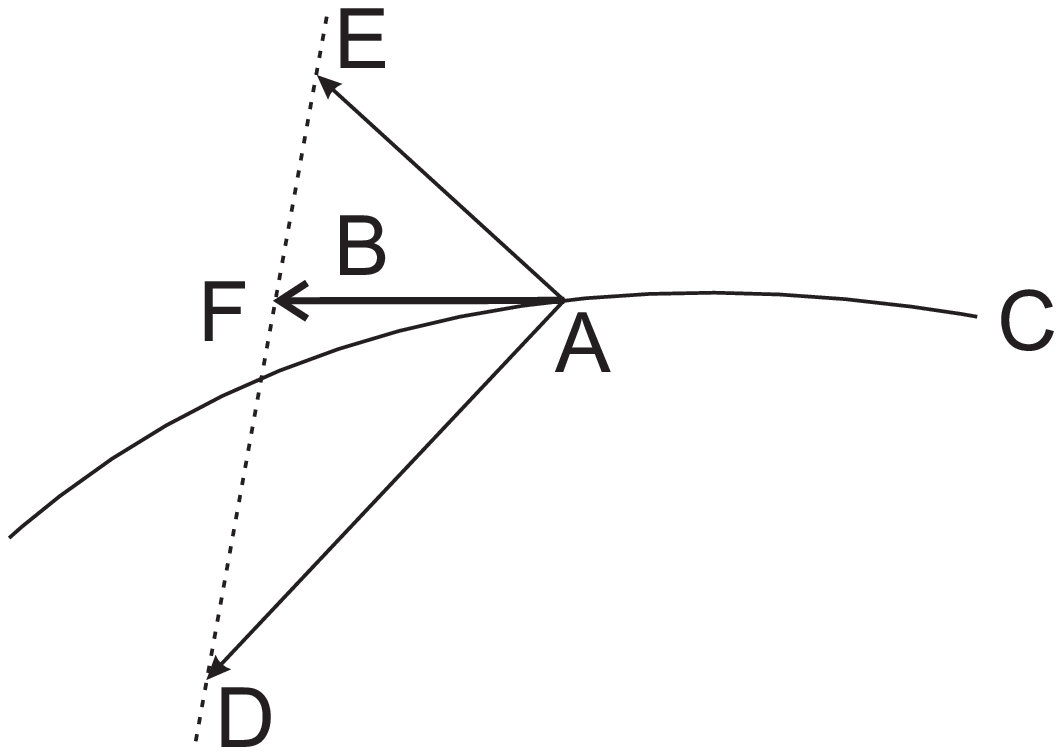}
\caption{\small{Fillipov's convention.}} \label{fig def filipov}
\end{center}
\end{minipage} \hfill
\begin{minipage}[b]{0.46\linewidth}
\begin{center}
\epsfxsize=5cm  \epsfbox{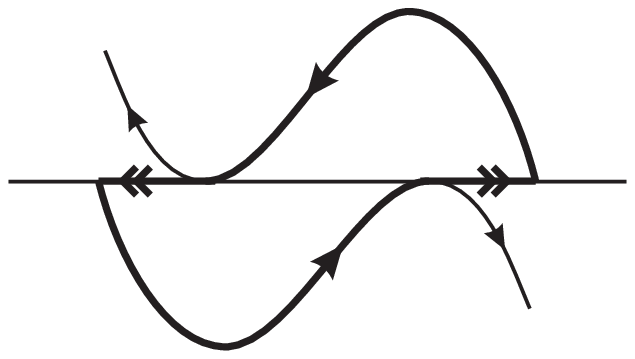}  \caption{\small{Canard
cycle.}}\label{fig canard}\end{center}
\end{minipage}
\end{figure}

Consider $X_0 \in \Omega^r.$ The \textit{sliding vector field}
associated to $X_0$ is the vector field  $X_0^s$ tangent to
$\Sigma_3$ and defined at $q\in \Sigma_3$ by $X_0^s(q)=m-q$ with $m$
being the point where the segment joining $q+X_1(q)$ and $q+X_2(q)$
is tangent to $\Sigma_3$ (see Figure \ref{fig def filipov}). It is
clear that if $q\in \Sigma_3$ then $q\in \Sigma_2$ for $-X_{0}$ and
then we  can define the {\it escaping vector field} on $\Sigma_2$
associated to $X_0$ by $X_0^e=-(-X_0)^s$. In what follows
we use the notation $X_0^\Sigma$ for both cases.\\

Our main interest here is to study a special kind of typical minimal
sets of non-smooth vector fields which will be called non-smooth
``\emph{canard cycles}" (see Figure \ref{fig canard}). A canard
cycle is a graphic composed by pieces of orbit of $X_1$, pieces of
orbit of the sliding vector field $X_0^\Sigma$ and/or pieces of
orbit of $X_2$. See Section 2 for a more precise definition.\\

An approximation of the non-smooth vector field $X_0 = (X_1, X_2 )$
by a $1$-parameter family $X_{\epsilon}$ of smooth vector fields is
called an $\epsilon$-regularization of $X_0$. We give the details
about this process in section 4. A transition function is used to
average $X_1$ and $X_2$ in order to get a family of smooth vector
fields that approximates $X_0$. The main aim is to deduce certain
dynamical properties of the non-smooth dynamical system from the
regularized system. What is familiar may or may not be a matter of
taste, at least it depends a lot on the dynamical properties of
one's interest. The regularization process developed by Sotomayor
and Teixeira produces a singular problem for which the discontinuous
set is a center manifold. Via a blow up we establish a bridge
between non-smooth systems and the geometric singular perturbation
theory.

Roughly speaking, the main results of this paper are the following:

\begin{itemize}
\item In our first result (Theorem \ref{teoA}), for a subclass of
non-smooth vector fields, we provide necessary and sufficient
conditions for the existence of canard kind solutions.

\item In our second result (Theorem \ref{teoB}),
following the ideas exposed in \cite{SM}, we prove that hyperbolic
canard cycles are limit sets, according Hausdorf distance, of
families of (smooth) hyperbolic limit cycles (this fact is not
proved in \cite{SM}). The regularization process plus a blow up
produce a singular perturbation problem $\mathcal{P_\epsilon}$. Our
result implies that the canard cycle is the periodic limit set of
closed orbits of $\mathcal{P_\epsilon}$, with $\epsilon \rightarrow
0$. An open problem is to use  the geometric singular perturbation
theory proposed by Dumortier and Roussarie
 (center manifolds obtained via saturation by the flow plus
blow up techniques, see \cite{DR} for details) to obtain the same
result.

\item In our third result (Theorem \ref{teoC}) we found an analogous
for {\it Poincar\'{e} Index} in the case of non-smooth vector
fields.
\end{itemize}

%

\section{Preliminaries and statements of the main results}

Consider $X_0 \in \Omega^r.$ We say that $q\in\Sigma$ is a
\textit{$\Sigma$-regular point} if
\begin{itemize}
\item [(i)] $X_1.f(q)X_2.f(q)>0$
or
\item [(ii)] $X_1.f(q)X_2.f(q)<0$ and $X_{0}^{\Sigma}(q)\neq0$ (that is $q\in\Sigma_2\bigcup\Sigma_3$ and it is not a singular
point of $X_{0}^{\Sigma}$).\end{itemize}

The points of $\Sigma$ which are not $\Sigma$-regular are called
\textit{$\Sigma$-singular}. We distinguish two subsets in the set of
$\Sigma$-singular points: $\Sigma^c$ and $\Sigma^f$. We say that $q
\in \Sigma^f$ is a \textit{pseudo equilibrium of $X_{0}$} if
$X_{0}^{\Sigma}(q)=0$ and we say that $q \in \Sigma^c$ is a
\textit{$\Sigma$-contact point} if $X_{0}^{\Sigma}(q) \neq 0$ and
$X_1.f(q)X_2.f(q) =0$ ($q$ is a contact point of $X_{0}^{\Sigma}$).\\

A $\Sigma$-contact point $q\in\Sigma^c$ is a \textit{$\Sigma$-fold
point} of $X_1$ if $X_1.f(q)=0$ but $X_1^{2}f(q)\neq0.$ Moreover,
$q\in\Sigma$ is a \textit{visible} (resp. {\it invisible})
\textit{$\Sigma$-fold point} of $X_1$
 if $X_1.f(q)=0$ and $X_1^{2}.f(q)> 0$
(resp. $X_1^{2}.f(q)< 0$). We say that $q$  is a $\Sigma$-fold point
 of
$X_{0}$ if it is a $\Sigma$-fold point either of $X_1$ or of $X_2$. \\

A pseudo equilibrium $q \in \Sigma^f$ is a \textit{$\Sigma$-saddle}
provided one of the following condition is satisfied: (i)
$q\in\Sigma_2$ and $q$ is an attractor for $X_{0}^{\Sigma}$ or (ii)
$q\in\Sigma_3$ and $q$ is a repeller for $X_{0}^{\Sigma}$. A pseudo
equilibrium $q\in\Sigma$ of $X_{0}$ is a $\Sigma$-\textit{repeller}
(resp. $\Sigma$-\textit{attractor}) provided $q\in\Sigma_2$ (resp.
$q \in \Sigma_3$) and $q$ is a repeller (resp. attractor) for
$X_{0}^{\Sigma}$. A point $q\in\Sigma$ is a \textbf{hyperbolic
pseudo equilibrium}
of $X_{0}$ if $q$ is a hyperbolic equilibrium point of $X_{0}^{\Sigma}.$\\

\begin{figure}[!h]
\begin{minipage}[b]{0.45\linewidth}
\begin{center}
\epsfxsize=4cm \epsfbox{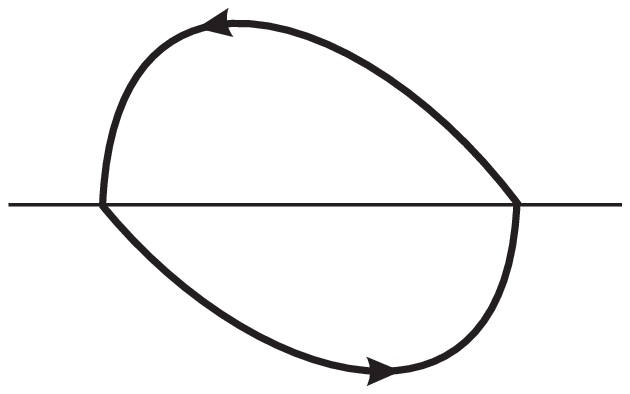} \caption{\small{Canard
cycle of kind I.}} \label{fig canard I}
\end{center}
\end{minipage} \hfill
\begin{minipage}[b]{0.45\linewidth}
\begin{center}
\psfrag{A}{$\Sigma=\Gamma$}\epsfxsize=3cm \epsfbox{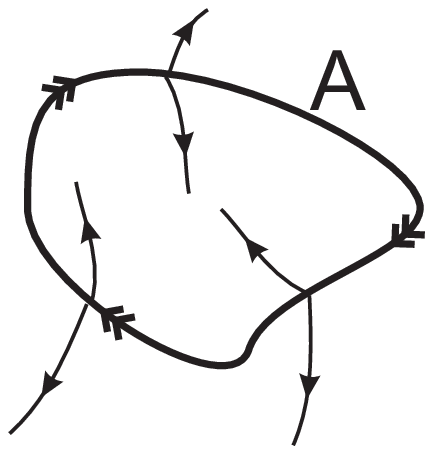}
\caption{\small{Canard cycle of kind II.}}\label{fig canard II}
\end{center}
\end{minipage}
\end{figure}

\begin{definition} Consider $X_{0} \in \Omega^r.$
\begin{enumerate}
\item A curve $\Gamma$ is a \textbf{canard cycle} if
$\Gamma$ is closed and 

  \begin{itemize}
  \item $\Gamma$ contains arcs of at least two of the vector fields $X_1 |_{\Sigma_{+}}$, $X_2 |_{\Sigma_{-}}$ and $X_{0}^{\Sigma}$ or is composed by a single arc of $X_{0}^{\Sigma}$;

  \item the transition between arcs of $X_1$
  and arcs of $X_2$ happens in sewing points (and vice versa);

  \item the transition between arcs of $X_1$
  (or $X_2$) and arcs of $X_{0}^{\Sigma}$ happens through
  $\Sigma$-fold points or regular points in the escape or sliding arc, respecting the orientation. Moreover
if $\Gamma\neq\Sigma$ then there exists at least one visible
$\Sigma$-fold point on each connected component of
$\Gamma\cap\Sigma$.
  \end{itemize}

\item Let $\Gamma$ be a canard cycle of $X_{0}$. We say that

  \begin{itemize}
  \item $\Gamma$ is a \textbf{canard cycle of kind
  I} if $\Gamma$ meets $\Sigma$ just in sewing points;

  \item $\Gamma$ is a \textbf{canard cycle of kind
  II} if $\Gamma = \Sigma$;

  \item $\Gamma$ is a \textbf{canard cycle of kind
  III} if $\Gamma$ contains at least one visible $\Sigma$-fold point of $X_{0}$.
  \end{itemize}

 In Figures \ref{fig canard I}, \ref{fig
canard II} and \ref{fig canard} appear  canard cycles of kind I, II
and III respectively.

\item Let $\Gamma$ be  a canard cycle. We say that
$\Gamma$ is \textbf{hyperbolic} if

  \begin{itemize}
  \item $\Gamma$ is of kind I and $\eta'(p) \neq 1$
  where $\eta$ is the first return map defined on a segment $T$ with $p\in T\pitchfork\gamma$;

  \item $\Gamma$ is of kind II;

  \item $\Gamma$ is of kind III and or $\Gamma\cap\Sigma\subseteq\Sigma_1\cup\Sigma_2$ or $\Gamma\cap\Sigma\subseteq\Sigma_1\cup\Sigma_3$.
  \end{itemize}

\end{enumerate}
\end{definition}

In \cite{SM} is proved that the $\epsilon$-regularization of
non-smooth vector fields $X_{0}$ with hyperbolic canard cycles has
hyperbolic limit cycles.

\begin{definition}\label{definicao politraj foco e grafico}
Let $(\overrightarrow{A \, B})_{X_{i}}$ be an arc of $X_i$ joining
the visible $\Sigma$-fold point $A$ to the point $B= X_i \pitchfork
\Sigma$. We say that $(\overrightarrow{A \, B})_{X_{i}}$ has
\textbf{focal kind} if there is not $\Sigma$-fold points between $A$
and $B$ (see Figure \ref{fig componente foco simples}) and we say
that $(\overrightarrow{A \, B})_{X_{i}}$ has \textbf{graphic kind}
if it has only one $\Sigma$-fold point between $A$ and $B$ (see
Figure \ref{fig componente grafico}), $i=1,2$.

\begin{figure}[!h]
\begin{minipage}[!h]{0.475\linewidth}
\begin{center}
\psfrag{A}{$A$} \psfrag{B}{$B$}
 \epsfxsize=4.7cm \epsfbox{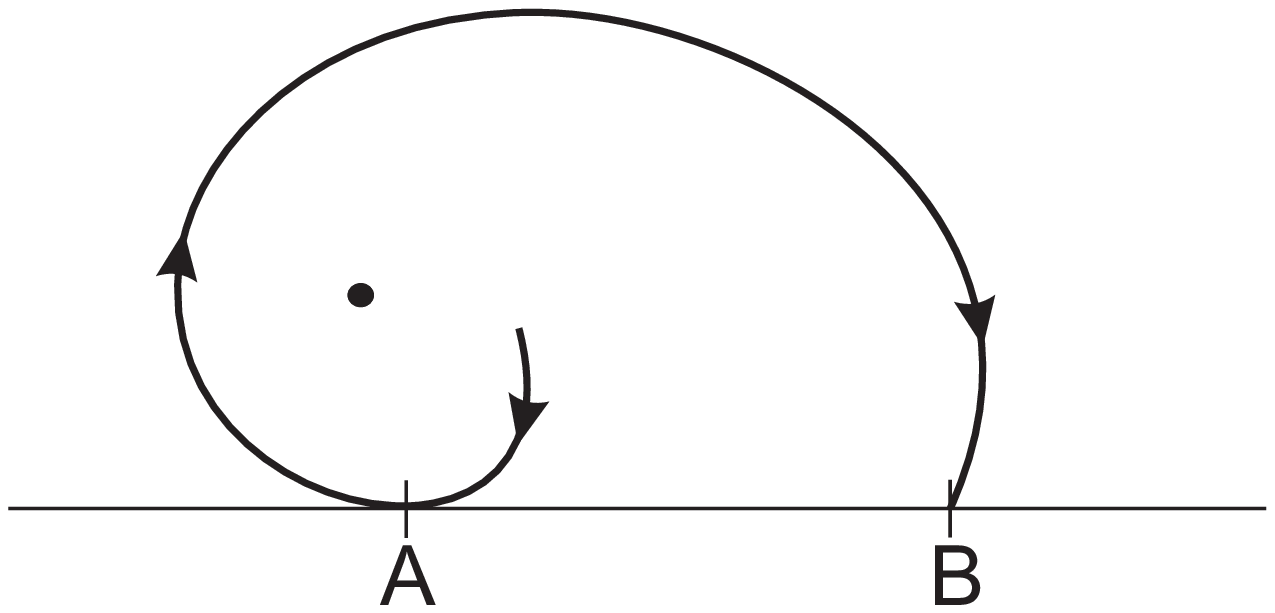}
\caption{\Small{Focal kind arc.}} \label{fig componente foco
simples}
\end{center}
\end{minipage} \hfill
\begin{minipage}[!h]{0.5\linewidth}
\epsfxsize=5.6cm \psfrag{A}{$A$} \psfrag{B}{$B$}
\epsfbox{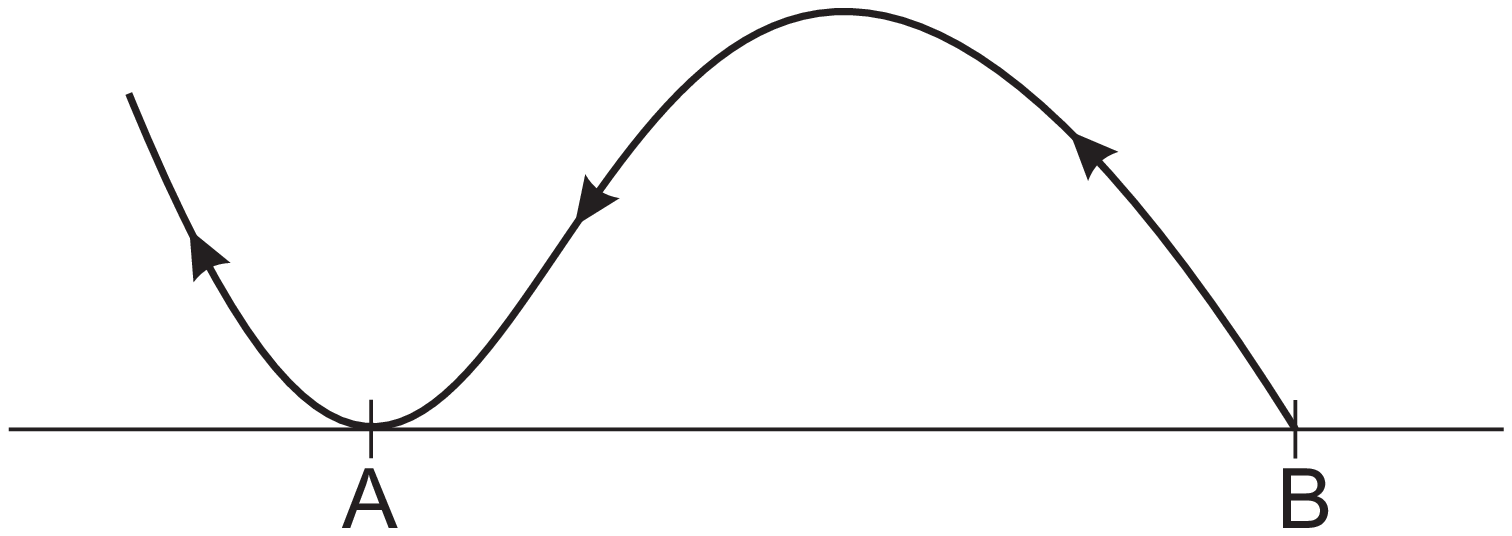} \caption{\Small{Graphic kind arc.}}
\label{fig componente grafico}
\end{minipage}
\end{figure}

\end{definition}

We remark that for a hyperbolic canard cycle we have that each
connected component of $\Gamma\cap\Sigma$ has only one $\Sigma$-fold
point (See \cite{SM} for more details).


By using the previous notation, our results are:

\begin{theorem}\label{teoA} Let
$X_{0}=(X_1,X_2)\in\Omega^r$ be a non-smooth vector field with $X_0$
presenting only one $\Sigma$-fold point $A$ which is visible. Denote
$\gamma_1$ the arc of $X_i$ ($i=1$ or $i=2$) which passes through
$A$ and  call $B$ the transversal contact point of $\gamma_{1}$ with
$\Sigma$. Then $X_{0}$ has a canard cycle $\Gamma$
 if and only if the following conditions are satisfied:
(i) the component $\gamma_{1}$ of $\Gamma$ which passes through $A$
is a focal kind arc; (ii) $X_{1}f. X_{2}f < 0$ in $(A,B]$ and (iii)
$\{ X_{1}, X_{2} \}$ is a linearly independent set in $[A,B]$.
Moreover, $\Gamma$ is of kind III.\end{theorem}

\begin{theorem}\label{teoB}Let $\Gamma_{0}$ be  a
hyperbolic canard cycle of $X_{0}$. Then for any $\epsilon>0$ the
regularized vector field $X_{\epsilon}$, has a hyperbolic limit
cycle  $\Gamma_{\epsilon}$ such that $ \Gamma_{\epsilon} \rightarrow
\Gamma_{0}$  when $\epsilon \rightarrow 0.$\end{theorem}

We remark that  the Hausdorff distance between compact sets of
$\R^2$ is: \vspace{-.15cm}\[D(K_1, K_2) = \max_{z_1 \in K_1,z_2 \in
K_2} \{d(z_1 , K_2) , d(z_2 , K_1)\}.\]

\begin{theorem}\label{teoC}Let $\Gamma_{0}$ be  a hyperbolic
canard cycle of the non-smooth vector field $X_{0}.$ If $\{p_{1}$,
$\ldots$, $p_{k}\}$ is the set of fixed or pseudo equilibrium points
(all hyperbolic) of $X_{0}$ inside $\Gamma_{0}$ then the index of
$\Gamma_{0}$ with respect to $X_{0}$ is the sum  of the
 index of $p_i$, for $i=1,...,k$. Moreover, this sum is equal to one.\end{theorem}

In section 5 we will define index of non-smooth vector fields.

The paper is organized as follows. In Sections 3, 4 and 5 we prove
Theorems \ref{teoA}, \ref{teoB} and \ref{teoC}, respectively. In
section 6 we apply Theorem \ref{teoA} to study a class of non-smooth
vector fields $X_0 \in \Omega^r$ with just one focal kind arc and
its bifurcation and we use the singular perturbation theory to study
hyperbolic canard cycles.


\section{Proof of the Theorem \ref{teoA}}

In this section we prove the first result of the paper.\\

\begin{figure}[!h]
\begin{center}
\psfrag{A}{} \psfrag{B}{$\gamma_{1}$} \psfrag{C}{$A$}
\psfrag{D}{$B$} \psfrag{E}{$X_{1}$} \psfrag{F}{$X_{2}$}
\psfrag{G}{$\Sigma$} \epsfxsize=5.8cm
\epsfbox{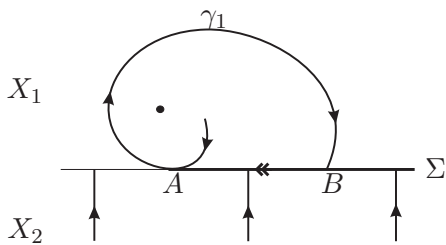} \caption{\small{Canard cycle with
just one $\Sigma$-fold point.}} \label{fig foco repulsor nova}
\end{center}
\end{figure}

\noindent\textbf{Proof.} First we prove that (i),(ii) and (iii)
imply the existence of the canard cycle. Since $X_{1}f . X_{2}f <0$
in $(A,B]$ the piece of $\Sigma$ between $A$ and $B$ is part of a
escaping region or a sliding region. Moreover since $\{ X_{1}, X_{2}
\}$ is a linearly independent set in $[A,B]$ the system does not
have pseudo equilibrium  points in $[A,B]$. Without lost of
generality, $[A,B]$ is part of the sliding region like in Figure
\ref{fig foco repulsor nova}. The curve $\Gamma = \gamma_{1} \cup
[B,A]$ is a hyperbolic canard cycle of kind III. We remark that
 this canard cycle  takes place in just one side of $\Sigma$.

 Now we prove that (i),(ii) and (iii) are necessaries conditions for
the existence of this particular kind of canard cycle. Since
$\Gamma$ is a hyperbolic canard cycle of kind III with just one
$\Sigma$-fold point, $\Gamma$ takes place in just one side of
$\Sigma$. In fact, if it does not occur, then $\Gamma$ returns to
$\Sigma$ at least twice and so  there exists at least a second
$\Sigma$-fold point. Without lost of generality we suppose that
$\Gamma$ is on the side corresponding to $X_1$. We denote by
$\gamma_1$ the part of the cycle $\Gamma$ which is a trajectory of
$X_1$. Thus we have that $\gamma_{1}$ is a focal  kind  arc because
if it is a graphic kind arc then there is another $\Sigma$-fold
point on $(A,B)$ (see Figure \ref{fig foco repulsor nova}). Since
$\Gamma$ has no one arc of $X_{2}$, the point $B$ belongs to an
escaping region or a sliding region and so $X_{1}f(B) . X_{2}f(B)
<0$. Let us assume that $B \in \Sigma_3$. Since $\gamma_1$ meets
$\Sigma$ in the point $B$, the flow slides via $X_0^\Sigma$ until
the point $A$ because there are not another $\Sigma$-fold point
between $A$ and $B$; therefore $X_{1}f . X_{2}f <0$ in $(A,B].$
Moreover, the linear independence of $\{ X_{1}, X_{2} \}$ on $[A,B]$
follows from the non-existence of pseudo equilibrium points on
$[A,B].$ \bbox
\\

Now, we will define an
auxiliar function which will be useful in the sequel.\\

Take $(A,B) \subset \Sigma_2 \cup \Sigma_3$ contained in the
escaping or in the sliding region. In $(A,B)$ consider  the point
$C=(C_{1},C_2)$, the vectors $X_1(C)=(D_1,D_2)$ and
$X_2(C)=(E_1,E_2)$ (as illustrated in Figure \ref{fig funcao
direcao}). The straight segment passing through $C+X_1(C)$ and $C +
X_2(C)$ meets $\Sigma$ in a point $p(C)$. We define the
C$^r$-application
$$
\begin{array}{cccc}
  p: & (A,B) & \longrightarrow & \Sigma \\
     & z & \longmapsto & p(z).
\end{array}
$$
We can choose local coordinates such that $\Sigma$ is the $x$-axis;
so $C=(C_1,0)$ and $p(C) \in \R \times \{ 0 \}$. The
\textit{direction function} on $\Sigma$ is defined by
$$
\begin{array}{cccc}
  H: & (A,B) & \longrightarrow & \R \\
     & z & \longmapsto & p(z) - z.
\end{array}
$$
\begin{figure}[!h]
\begin{center}
\psfrag{A}{$A$} \psfrag{B}{$B$} \psfrag{C}{$C$} \psfrag{D}{$\Sigma$}
\psfrag{E}{$X_1$} \psfrag{F}{$X_2$} \psfrag{G}{$C+X_2(C)$}
\psfrag{H}{$C+ X_1(C)$} \psfrag{I}{$p(C)$} \epsfxsize=5cm
\epsfbox{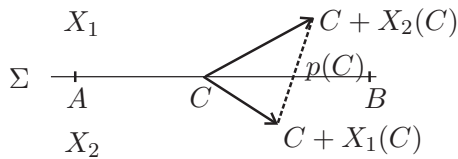} \caption{\small{Direction function.}}
\label{fig funcao direcao}
\end{center}
\end{figure}

We have that $H$ is a C$^r$-application and

\begin{itemize}
\item if $H(C) < 0$ then the orientation of $X_0^{\Sigma}$ in a small neighborhood of $C$ is from $B$ to $A$;

\item if $H(C) = 0$ then  $C \in \Sigma^f$;

\item if $H(C) > 0$ then the orientation of $X_0^{\Sigma}$ in a small neighborhood of $C$ is from $A$ to $B$.
\end{itemize}

Simple calculations shows that $p(C_1) = \frac{(D_1+C_1) (E_2) -
(D_2) (E_1+C_1)}{(E_2) - (D_2)}$.\\

Assuming all the hypothesis of Theorem \ref{teoA} we have the
following corollary.

\begin{corollary}\label{corteoA} The non-smooth vector field $X_{0}$ has a canard cycle $\Gamma$
 if and only if the direction function $H: [A,B]
\longrightarrow \R$ is a well defined function and it has no zeros.
Moreover, $\Gamma$ is of kind III.\end{corollary}


\section{Proof of Theorem \ref{teoB}}

First of all we present the concept of $\epsilon$-regularization of
non-smooth vector fields. It was introduced by Sotomayor and
Teixeira in \cite{ST}. The regularization  gives the mathematical
tool to study the stability of these systems, according to the
program introduced by Peixoto. The method consists in the analysis
of the regularized vector field which is  a smooth approximation of
the non-smooth vector field. Using this process we get a
$1$-parameter family of vector fields $X_{\e}\in \chi^r(K,\R^2)$
such that for each $\e_0 > 0$ fixed we have
\begin{itemize}
\item [(i)]   $X_{\e_0}$ is equal to $X_1$ in all
points of $ \Sigma_+$ whose  distance to $\Sigma$ is bigger than
$\e_0;$
\item [(ii)] $X_{\e_0}$ is equal to $X_2$ in all points of $ \Sigma_-$ whose distance to $\Sigma$ is bigger than $\e_0$.
\end{itemize}

\begin{definition}
 A $C^\infty$ function $\varphi:\R \longrightarrow \R$
is a transition function if $\varphi(x)=-1$ for $x\leqslant -1$,
$\varphi(x)=1$ for $x\geqslant 1$ and $\varphi'(x)>0$ if
$x\in(-1,1).$ The $\epsilon$-regularization of $X_0=(X_1,X_2)$ is
the 1-parameter family $X_{\e}\in \chi^r$ given by
\begin{equation} \label{regularization}
X_{\e}(q)=\left(
\dfrac{1}{2}+\dfrac{\varphi_{\e}(f(q))}{2}\right)X_1(q) +\left(
\dfrac{1}{2}-\dfrac{\varphi_{\e}(f(q))}{2}\right) X_2(q).
\end{equation}
with $\varphi_{\e}(x)=\varphi(x/\e),$ for $\e>0.$
\end{definition}

In order to prove Theorem \ref{teoB} we need to construct a special
neighborhood of arbitrary diameter for hyperbolic canard cycles.\\

\noindent\textbf{Construction of a neighborhood of diameter $\mu$
around a hyperbolic canard cycle.} Here we describe a method to
construct a tubular neighborhood of diameter $\mu$ around a
hyperbolic canard cycle. This presentation is done for canard cycles
of kind III, but the ideas can also be extended for kinds I or II.
We will be particularly interested in two of them: the ones that
take place on just  one side of $\Sigma$ and with just one visible
$\Sigma$-fold point and the ones that take place on the two sides of
$\Sigma$ with two visible $\Sigma$-fold points (one for $X_1$ and
another one for $X_2$).\\

\vspace{-.1cm}

\textit{Case 1- One $\Sigma$-fold point.} Denote by $\Gamma$ the
hyperbolic canard cycle of kind III with just one $\Sigma$-fold
point and with orientation showed in Figure \ref{viz tubular 1
dobra} (the reverse orientation is treated in a similar way).
Consider the strip of diameter $\mu$ around $\Sigma = \{ y=0  \}$.
Let  $p_{1}$ and $q_{1}$ be points in $\{ y=\mu \} \cap \Gamma$.
Take an arc $\gamma_2$ of the vector field $X_1$ passing to the
point $p_{2} \in \{ y=\mu \}$ in such a way that $p_2$ stays on the
left of $p_1$ and such that $\gamma_{2}$ returns to the line $ y=\mu
$ in a point $q_{2}$ which is in a neighborhood of $q_{1}$. Take
this trajectory satisfying $d( \Gamma, \gamma_{2} ) < \frac{\mu}{2}$
(this is possible by the continuity of $X_1$). Let $r_{1}$ be the
point where the arc of $X_1$  through by $p_{1}$ first meets the
straight line $y = \mu$ for negative time. Analogously  take an arc
$\gamma_1$ of the field $X_1$ passing by the point $r_{2}$ in such a
way that $r_2$ stays on the left of $r_1$ and such that $\gamma_{1}$
has second return to $y=\mu$ in a point $q_{3}$ which is in a
neighborhood of $q_{1}$. Take this trajectory satisfying $d( \Gamma,
\gamma_{1} ) < \frac{\mu}{2}$. On $\Sigma$, on the left of the
$\Sigma$-fold point $A$, the flow of $X_1$ is oriented to up, so it
is possible to construct a transversal section $\sigma_{2}$ joining
$p_{2}$ to the straight line $y=0$ in such a way that the same
trajectories of $X_1$  cross transversally $\sigma_{2}$ and the
segment $\overline{p_{2} \, p_{1}}$. Take $t_{2}$  the point where
$\sigma_{2}$ meets the straight line $y=0$  satisfying $d(\Gamma,
\sigma_{2})< \mu$, as before. Moreover, on $\Sigma$, on the right of
the point $B$, the flow of $X_{1}$ is oriented for down, so it is
possible to construct a transversal section $\theta_{2}$ joining
$q_{2}$ to the straight line $y=0$ in such a way that the
trajectories of $X_1$ that cross transversally $\theta_{2}$ do not
cross the segment $\overline{q_{1} \, q_{2}}$. Let  $s_{2}$ be the
point where $\theta_{2}$ meets the straight line $y=0$ (here we also
need to take care for $d(\Gamma, \theta_{2})< \mu$). Since $[A, B]$
is a sliding region, the flow of $X_{2}$ is transversal to $\Sigma$.
In the straight line $y=-\mu$, consider the points $u_{2}$ and
$v_{2}$ ($u_{2}$ is on the left of $v_{2}$) satisfying that the
trajectories of $X_{2}$ crossing the transversal sections
$\lambda_{2} = \overline{t_{2} \, u_{2}}$ and $\delta_{2} =
\overline{s_{2} \, v_{2}}$ meet transversally $\Sigma$ at the
segment $\overline{t_{2} \, s_{2}}$ (again, we need to take care for
$d(\Gamma , \lambda_{2}) < \mu$ and $d(\Gamma , \delta_{2})<\mu$).

\begin{figure}[!h]
\begin{center}
\psfrag{A}{$A$} \psfrag{B}{$B$} \psfrag{C}{$p_{1}$}
\psfrag{D}{$p_{2}$} \psfrag{E}{$\gamma_{2}$}
\psfrag{F}{$\gamma_{1}$} \psfrag{G}{$\Gamma$} \psfrag{H}{$r_{2}$}
\psfrag{I}{$r_{1}$} \psfrag{J}{$\theta_{2}$} \psfrag{K}{$q_{1}$}
\psfrag{L}{$q_{3}$} \psfrag{M}{$q_{2}$} \psfrag{N}{$\delta_{2}$}
\psfrag{O}{$\Sigma$} \psfrag{P}{$\mu$} \psfrag{Q}{$\lambda_{2}$}
\psfrag{R}{$\sigma_{2}$} \psfrag{S}{$s_{2}$} \psfrag{T}{$t_{2}$}
\psfrag{V}{$v_{2}$} \psfrag{U}{$u_{2}$} \psfrag{X}{$X_{1}$}
\psfrag{Z}{$X_{2}$} \epsfxsize=8cm \epsfbox{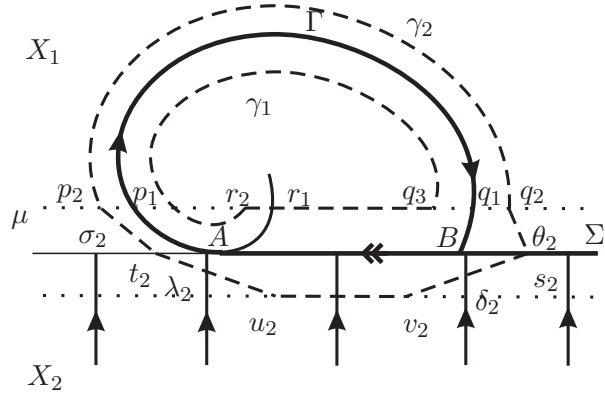}
\caption{\small{Tubular neighborhood of a canard cycle with one
$\Sigma$-fold.}} \label{viz tubular 1 dobra}
\end{center}
\end{figure}

In this way, the strip defined by the closed curve $\gamma_{1} \cup
\overline{q_{3} \, r_{2}}$ and by the closed curve $\gamma_{2} \cup
\theta_{2} \cup \delta_{2} \cup \overline{v_{2} \, u_{2}} \cup
\lambda_{2} \cup \sigma_{2}$ is a tubular neighborhood of $\Gamma$
of diameter $\mu$. Note that the flow of $X_{0}=(X_{1},X_{2})$ is
arriving in this neighborhood and never it departs from it.\\

\begin{figure}[!h]
\begin{center}
\psfrag{A}{$\mu$} \psfrag{B}{$A'$} \psfrag{C}{$\rho_{1}$}
\psfrag{D}{$\Sigma$} \psfrag{E}{$e$} \psfrag{F}{$f$}
\psfrag{G}{$\Gamma$} \psfrag{H}{$q_{1}$} \psfrag{I}{$i$}
\psfrag{J}{$j$} \psfrag{K}{$u_{1}$} \psfrag{L}{$u_{2}$}
\psfrag{M}{$\eta_{2}$} \psfrag{N}{$\eta_{1}$} \psfrag{O}{$q_{2}$}
\psfrag{P}{$\pi_{1}$} \psfrag{Q}{$\gamma_{1}$} \psfrag{R}{$r_{2}$}
\psfrag{S}{$\sigma_{1}$} \psfrag{T}{$t_{1}$} \psfrag{U}{$\pi_{2}$}
\psfrag{V}{$v_{1}$} \psfrag{Y}{$p_{1}$} \psfrag{X}{$x_{2}$}
\psfrag{Z}{} \psfrag{1}{$h$} \psfrag{2}{$\gamma_{2}$}
\psfrag{3}{$\sigma_{2}$} \psfrag{4}{$\rho_{2}$} \psfrag{5}{$A$}
\psfrag{6}{$t_{2}$} \psfrag{7}{$X_{2}$} \psfrag{8}{$p_{2}$}
\psfrag{9}{$v_{2}$} \psfrag{10}{$\theta_{2}$} \psfrag{11}{}
\psfrag{12}{$x_{1}$} \psfrag{13}{$r_{1}$} \psfrag{14}{$\theta_{1}$}
\psfrag{15}{$c$} \psfrag{16}{$g$} \psfrag{17}{$X_{1}$}
\psfrag{0}{$d$} \psfrag{}{} \psfrag{}{} \psfrag{}{} \psfrag{}{}
\psfrag{}{} \psfrag{}{}
 \psfrag{}{} \psfrag{}{} \psfrag{}{} \psfrag{}{}
\epsfxsize=10cm \epsfbox{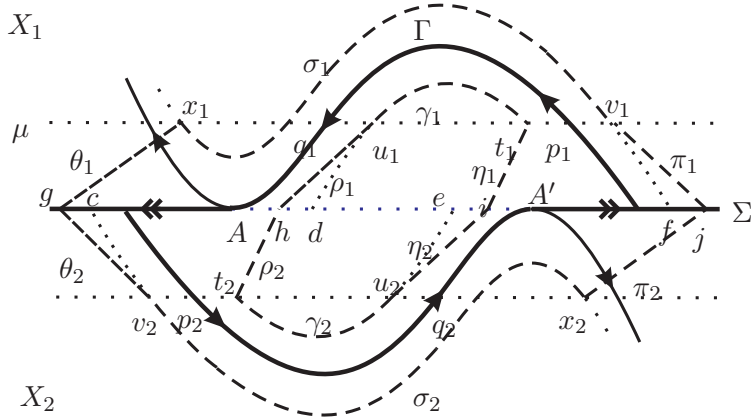}
\caption{\small{Tubular neighborhood for canard cycles of kind III
with two $\Sigma$-fold points.}} \label{viz tubular 2 dobra}
\end{center}
\end{figure}

\textit{Case 2- Two $\Sigma$-fold points.} Now we study the
hyperbolic canard cycles of kind III with two visible $\Sigma$-fold
points, being one for $X_1$ and the other one for $X_2$, like showed
in Figure \ref{viz tubular 2 dobra}. We work with canard cycles
$\Gamma$ that have only escaping regions on  $\Sigma = \{ y=0 \}$
(the case with sliding regions is treated similarly). Consider the
strip of diameter $\mu$ around $\Sigma$. Let $p_{1}$ and $q_{1}$ be
points in $\{ y=\mu \} \cap \Gamma$. Take an arc $\gamma_1$ of the
vector field $X_1$ through $t_{1} \in \{ y=\mu \}$ satisfying that
$t_1$ stays on the left of $p_1$ and such that $\gamma_{1}$ returns
to the line $ y=\mu $ in a point $u_{1}$ which is in a neighborhood
of $q_{1}$. Take this trajectory satisfying that $d( \Gamma,
\gamma_{1} ) < \frac{\mu}{2}$. Take an arc $\sigma_1$ of the vector
field $X_1$ through $v_{1} \in \{ y=\mu \}$ satisfying that $v_1$
stays on the right of $p_1$ and such that  $\sigma_{1}$ has second
return on the straight line $y=\mu$ in a point $x_{1}$, even take
this trajectory with the particularity that $d( \Gamma, \sigma_{1} )
< \frac{\mu}{2}$. We repeat the same argument for the vector field
$X_{2}$ and we found the points $p_{2}$, $q_{2}$, $t_{2}$, $u_{2}$,
$v_{2}$ and $x_{2}$ respectively, and the curves $\gamma_{2}$ and
$\sigma_{2}$. Let  $c$ be the point on $\Sigma\cap\sigma_2$, $d$ be
the point on $\Sigma\cap\gamma_1$, $e$ be the point on
$\Sigma\cap\gamma_2$ and $f$ be the point on $\Sigma\cap\sigma_1$ as
indicated in Figure \ref{viz tubular 2 dobra}. On $\Sigma$, take the
points $g$ on the left of $c$, $h$ between $A$  and $d$, $i$ between
$e$ and $A'$  and $j$ on the right of $f$; satisfying that the arcs
$\theta_{1}$ (joining $g$ to $x_{1}$), $\rho_{1}$ (joining $u_{1}$
to $h$), $\eta_{1}$ (joining $i$ to $t_{1}$) and $\pi_{1}$ (joining
$v_{1}$ to $j$) are transversal sections for $X_{1}$ and the arcs
$\theta_{2}$ (joining $g$ to $v_{2}$), $\rho_{2}$ (joining $h$ to
$t_{2}$), $\eta_{2}$ (joining $u_{2}$ to $i$) and $\pi_{2}$ (joining
$j$ and $x_{2}$) are transversal sections for $X_{0}$ with the
distance from $\Gamma$ to any one of this arcs less than $\mu$.

In this way, the strip defined by the closed curve $\gamma_{1} \cup
\rho_{1} \cup \rho_{2} \cup \gamma_{2} \cup \eta_{2} \cup \eta_{1}$
and by the closed curve $\sigma_{1} \cup \theta_{1} \cup \theta_{2}
\cup \sigma_{2} \cup \pi_{2} \cup \pi_{1}$ is a tubular neighborhood
for $\Gamma$ of diameter $\mu$. Note that the flow of
$X_{0}=(X_{1},X_{2})$ is departing from the tubular neighborhood and
never it arrives in it.

\begin{itemize}

\item Since this neighborhood bounds a region where the non-smooth vector field
$X_{0} = (X_{1} , X_{2})$ is arriving in or it is departing from
them, it makes sense to say \textit{attractor canard cycle} or
\textit{repeller canard cycle}.

\item In the neighborhoods constructed before we allow that
trajectories can make part of them, however it is possible to do it
with the  flow of $X_0$ being transversal to the boundaries of the
tubular neighborhoods. In fact, it is enough to replace the
trajectories by transversal curves. It is important for the
construction of the tubular neighborhood of the canard cycles of
kind I. Thus we make a construction like we made before but now we
can use for this, the first return application $\eta$ and thus if
$\eta'<1$ we have an attractor canard cycle and if $\eta'>1$ we have
a repeller canard cycle.

\item For canard cycles of kind II is enough to take the strip
of diameter $\mu$ in the beginning of the construction as the tubular neighborhood.

\item Any other hyperbolic canard cycle is  an arrangement of
pieces of the canard cycles described above and so we can construct
a tubular neighborhood for it arranging the previous tubular
neighborhoods.\\
\end{itemize}

\noindent\textbf{Proof of Theorem \ref{teoB}}. Let $\Gamma_{0}$ be a
canard cycle of $X_{0}$ and let $V_\e$ be a tubular neighborhood of
diameter $\e$ around $\Gamma_0$. Since $X_0$ is transversal to the
boundary of $V_\e$, by continuity, the regularized vector field
$X_\e$ also is transversal to the boundary of $V_\e$. Assume that
$\Gamma_0$ is an attractor canard cycle, so the flow of $X_0$ is
arriving in the neighborhood $V_\e$ and consequently the flow of
$X_\e$ also is arriving in the neighborhood $V_\e$. As there are not
fixed points in $V_\e$, applying the Poincar\'{e}-Bendixson Theorem
we conclude that there exists an attractor limit cycle $\Gamma_\e$
inside $V_\e.$ Moreover with a more detailed analysis we can prove
that it is hyperbolic (see \cite{SM} for instance). Since every
paths that compose $V_\e$ depends continuously of $\e$ we have that
the diameter of the tubular neighborhood is a continuous function of
the variable $\e$. Therefore making $\epsilon \rightarrow 0$ we
conclude that $\Gamma_{\epsilon} \rightarrow \Gamma_{0}$ (see Figure
\ref{conv ciclos
2}).\bbox\\

\vspace{-.5cm}\begin{figure}[!h]
\begin{center}
\psfrag{A}{} \psfrag{B}{} \psfrag{C}{} \psfrag{D}{$\Sigma$}
\psfrag{E}{} \psfrag{F}{$\epsilon$} \psfrag{G}{$\Gamma_{0}$}
\psfrag{H}{$\Gamma_{\epsilon}$} \psfrag{I}{} \psfrag{V}{$V_{0}$}
\psfrag{U}{$V_{\epsilon}$ } \psfrag{}{} \epsfxsize=4.5cm
\epsfbox{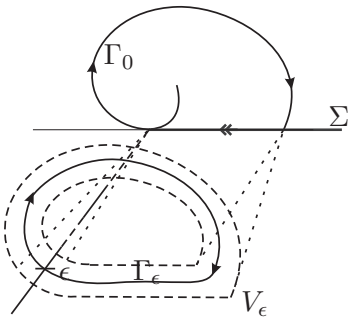} \caption{\small{Cycles
convergence.}} \label{conv ciclos 2}
\end{center}
\end{figure}

We remark that if $\Gamma_{0}$ is an attractor (resp. repeller)
hyperbolic canard cycle of $X_{0}$, then the same occurs for
$\Gamma_{\epsilon}$ and $X_{\epsilon}.$

\section{Proof of Theorem \ref{teoC}}Now we start our discussion about the third result.
 Let $I=[0,1]$ be an interval
and $\sigma: I \rightarrow \mathcal{U}$ be an oriented closed
continuous path. Suppose that there are no critical points of $X_0$
on $\sigma$. Let us move a point $P$ along the curve in the
counterclockwise direction. The vector $X_0(P)$ will rotate during
the motion. When $P$ returns to its starting place after one
revolution along the curve $\sigma$, $X_0(P)$ also returns to its
original position. During the journey $X_0(P)$ will make some whole
number of revolutions. Counting these revolutions positively if they
are counterclockwise, negatively if they are clockwise, the
resulting algebraic sum of the number of revolutions is called the
\textit{index of $\sigma$ with respect to $X_0$}, and is denoted by
$I(X_0,\sigma)$.

To calculate $I(X_0, \sigma)$ is convenient normalize $X(\sigma(t))$
as an unit vector at the origin. In this way, we can define a
function $\theta: I \rightarrow \R$ such that
\[\displaystyle\lim_{\overline{t}
\rightarrow t^{-}} \frac{X(\sigma(\overline{t}))}{\|
X(\sigma(\overline{t})) \|} = \displaystyle\lim_{\overline{t}
\rightarrow t^{-}} (\cos \theta(\overline{t}), \sin
\theta(\overline{t}))\] for every $\overline{t} \in I$. The function
$\theta$ is called \textit{angle function}.

We observe that in the case of smooth vector fields, the angle
function is always continuous, but in the case of non-smooth vector
field it admits a ``jump" when the path pass to a point $s_i \in
\Sigma$, $i \in \N$. Therefore, we establish a rule for this jump;
at $s_i = \sigma(t_i)$ the angle function oscillates from
$\displaystyle\lim_{t \rightarrow t_i^{-}}\theta(t)$ to
$\displaystyle\lim_{t \rightarrow t_i^{+}}\theta(t)$. If
$\displaystyle\lim_{t \rightarrow t_i^{-}}\theta(t)\in
I_{k-1}=(2(k-1)\pi, 2k \pi)$ and $\displaystyle\lim_{t \rightarrow
t_i^{+}}\theta(t)\in I_k =(2k \pi, 2(k+1) \pi)$, where $k \in \Z$,
we add $1$ to the the number $I(X_0, \sigma)$; if
$\displaystyle\lim_{t \rightarrow t_i^{-}}\theta(t)\in I_{k}$ and
$\displaystyle\lim_{t \rightarrow t_i^{+}}\theta(t)\in I_{k-1}$  we
add $-1$ to the number $I(X_0, \sigma)$. We always consider that the
jump of the vector $\displaystyle\lim_{t \rightarrow
t_i^{-}}X(\sigma(t))$ to the vector $\displaystyle\lim_{t
\rightarrow t_i^{+}}X(\sigma(t))$ occurs by the smallest angle
between this vectors.


\begin{figure}[!h]
\begin{center}
\psfrag{A}{$X_2(s_1)$} \psfrag{B}{$\alpha_1$} \psfrag{J}{$s_1$}
\psfrag{C}{$X_1(s_2)$} \psfrag{D}{$\alpha_2$} \psfrag{E}{$\Sigma$}
\psfrag{F}{$\sigma(I)$} \psfrag{G}{$X_1(s_1)$} \psfrag{H}{$s_2$}
\psfrag{K}{$X_0^\Sigma(s_2)$} \psfrag{I}{$s_3$} \psfrag{U}{}
\psfrag{}{} \epsfxsize=7cm
\epsfbox{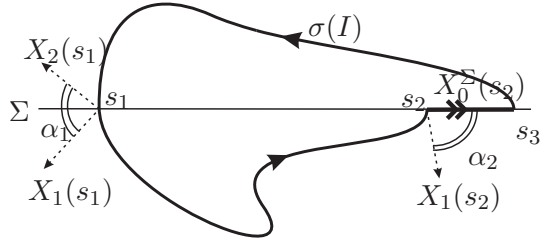} \caption{\small{Angle
function: at $s_1$ it has a jump of size $\alpha_1$ and at $s_2$ it
has a jump of size $\alpha_2$.}} \label{fig funcao angulo
descontinua}
\end{center}
\end{figure}

The difference $\theta(1) - \theta(0)$ is a multiple of $2 \pi$, and
$$
I(X_0, \sigma)  =  \frac{\theta(1) - \theta(0)}{2 \pi}
$$
is  an integer independent of the chosen $\sigma$-parametrization.
This number also is called \textit{Poincar\'{e} Index of the curve
$\sigma$ with relation to the non-smooth vector field $X_0$}.\\

\vspace{-.15cm}

Our interest here is to calculate the index of canard cycles
surrounding fixed or pseudo equilibrium points that are the critical
points of $X_0$. We will see that, different from the smooth case,
given two canard cycles $\Gamma_1$ and $\Gamma_2$ surrounding the
same critical points of $X_0$
we have $I(X_0,\Gamma_1) \neq I(X_0,\Gamma_2)$ in general.\\

\begin{figure}[!h]
\begin{center}
\psfrag{A}{$X_{1}$} \psfrag{B}{$X_{2}$} \psfrag{C}{$p_{1}$}
\psfrag{D}{$\Sigma$} \psfrag{E}{$a$} \psfrag{F}{$f$} \psfrag{G}{$g$}
\psfrag{H}{$h$} \psfrag{I}{$i$} \psfrag{J}{$\gamma_{1}$}
\psfrag{K}{$\delta_{1}$} \psfrag{L}{$e$} \psfrag{M}{$p_{2}$}
\psfrag{N}{$\gamma_{2}$} \psfrag{O}{$d$} \psfrag{P}{$\delta_{2}$}
\psfrag{Q}{$c$} \psfrag{R}{$b$} \epsfxsize=7cm
\epsfbox{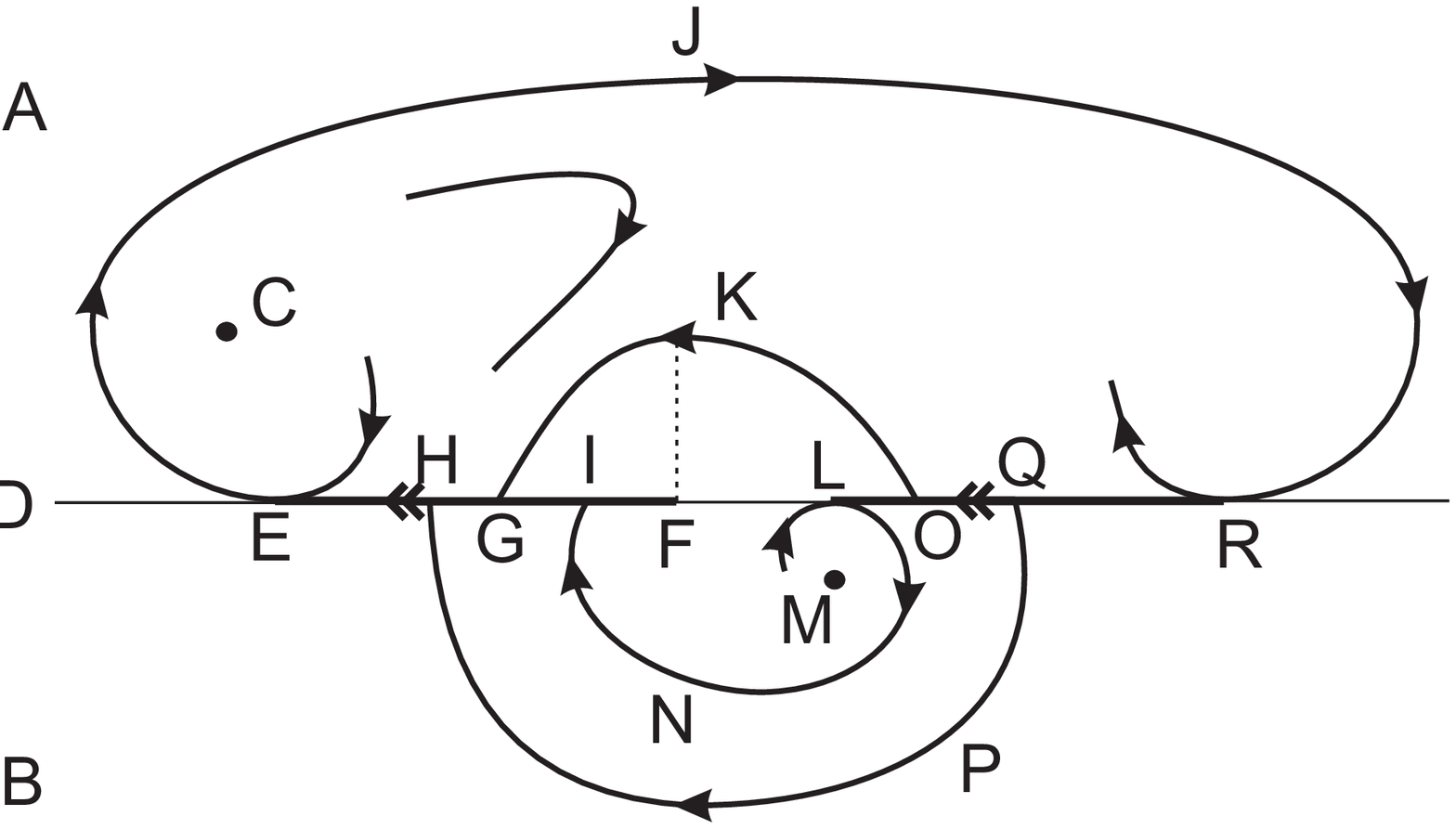} \caption{\small{Non-hyperbolic canard
cycle.}} \label{politrajetoria nao hiperbolica 1}
\end{center}
\end{figure}

\noindent\textbf{Example.} Consider the configuration described in
Figure \ref{politrajetoria nao hiperbolica 1}. Let $\gamma_{1}$ be
an arc of $X_{1}$ joining the $\Sigma$-fold points $a$ and $b$; $e$
the $\Sigma$-fold point of $X_2$; $c$ and $d$ points in the escaping
region $\overline{e \, b}$; $g$, $h$ and $i$ points in the sliding
region $\overline{f \, a}$ where $f$ is a invisible $\Sigma$-fold
point of $X_1$; $\delta_{2}$ the arc of $X_{2}$ joining $c$ and $h$,
$\delta_{1}$ the arc of $X_{1}$ joining $d$ and $g$; and
$\gamma_{2}$ the arc of $X_{2}$ joining $e$ and $i$. Consider the
fixed point of $X_{2}$ named by $p_{2}$ in Figure
\ref{politrajetoria nao hiperbolica 1}. We can choose closed curves
such that the Poincar\'{e} index to non-smooth vector fields is any
natural number. For example, take the path $\overline{h \, a} \cup
\gamma_{1} \cup \overline{b \, c} \cup \delta_{2}$ and the index is
$1$; in an analogous way if we take the path $\overline{i \, a} \cup
\gamma_{1} \cup \overline{b \, c} \cup \delta_{2} \cup \overline{h
\, a} \cup \gamma_{1} \cup \overline{b \, e} \cup \gamma_{2}$ the
index is $2$. Repeating this argument we can found closed curves
such that the Poincar\'{e} index
is any integer number that we wish.\\

The previous situation just occurs because the canard cycle given,
for example, by $\gamma_{1} \cup \overline{b \, c} \cup \delta_{2}
\cup \overline{h \, a}$, is non-hyperbolic because in its
composition we can found pieces of escaping region and pieces of
sliding region. If we eliminate this possibility we have the next
theorem, which has an analogous in the case of smooth vector fields.

\begin{remark}\label{obs convergencia singularidades} \emph{We
recall
that if $X_{0}$ has a hyperbolic $\Sigma$-saddle
 (or a hyperbolic $\Sigma$-focus) $s_{0}$ then the regularized vector field $X_{\epsilon}$ has a hyperbolic saddle (hyperbolic focus) $s_{\epsilon}$
 where
 $s_{\epsilon} \rightarrow s_{0}$ when
$\epsilon \rightarrow 0$ (for details see \cite{SM}). Moreover, we
can verify that if $s_0$ is a saddle or a $\Sigma$-saddle then we
can take a sufficiently small closed path $\sigma$ around $s_0$ and
prove that $I(X_0, \sigma)=-1$ (and if $s_0$ is a focus or a
$\Sigma$-attractor or a $\Sigma$-repeller then we can take a
sufficiently small closed path $\widetilde{\sigma}$ around $s_0$ and
prove that $I(X_0, \widetilde{\sigma})=1$). When the path $\sigma$
is sufficiently small to have just one critical point of $X_0$,
named $s_0$, in its interior we use the notation $I_{s_0}(X_0,
\sigma)$ to denote its Poincar\'{e} Index.}
\end{remark}

\subsection{Proof of Theorem \ref{teoC}}

We want to prove that if $\Gamma_{0}$ is a hyperbolic canard cycle
of the non-smooth vector field $X_{0}$ and if $p_{1}$, $\ldots$,
$p_{k}$ are the only ones critical points (all hyperbolic) of
$X_{0}$ inside $\Gamma_{0}$ then $I(X_0, \Gamma_0)$ is well defined
and
$$
I(X_0,\Gamma_0) = \sum_{i=1}^{k}I_{p_i}(X_0, \Gamma_0) = 1.
$$
First of all we assume that the index is well defined. Let $X_{0}$
be a non-smooth vector field with a hyperbolic canard cycle
$\Gamma_{0}$ and $k$ hyperbolic critical points of $X_0$ inside
$\Gamma_0$. Thus, the regularized vector field $X_{\epsilon}$ has a
hyperbolic limit cycle $\Gamma_{\epsilon}$ and $k$ hyperbolic fixed
points inside $\Gamma_\epsilon$. So, by the Poincar\'{e} Index
Theorem for smooth vector fields the index calculated in
$\Gamma_{\epsilon}$ in relation to $X_{\epsilon}$ is the sum of the
index of the fixed points of $X_{\epsilon}$ inside
$\Gamma_{\epsilon}$ and this sum is equal to $1$. Since
$\Gamma_{\epsilon} \rightarrow \Gamma_{0}$ we conclude that the
index calculated in $\Gamma_{0}$ in relation to $X_{0}$ is the sum
of the index of the critical points of $X_{0}$ inside $\Gamma_{0}$
and this sum is equal to 1 (see remark \ref{obs convergencia
singularidades}). In order to finish the proof we verify that the
index  is well defined in the case that the closed curve
$\Gamma_{0}$  is a hyperbolic canard cycle. In fact, let
$\Gamma_{0}$ be a hyperbolic canard cycle of $X_{0}$. Let us assume
that $\Gamma_{0}$ is of one kind described in Figures \ref{fig
politrajetoria foco atrator} or \ref{fig politrajetoria foco
repulsor} below.

\begin{figure}[!h]
\begin{minipage}[b]{0.485\linewidth}
\begin{center}
\psfrag{A}{$p_{1}$} \psfrag{B}{$\gamma_{1}$} \psfrag{C}{$A$}
\psfrag{D}{$B$} \epsfxsize=5cm \epsfbox{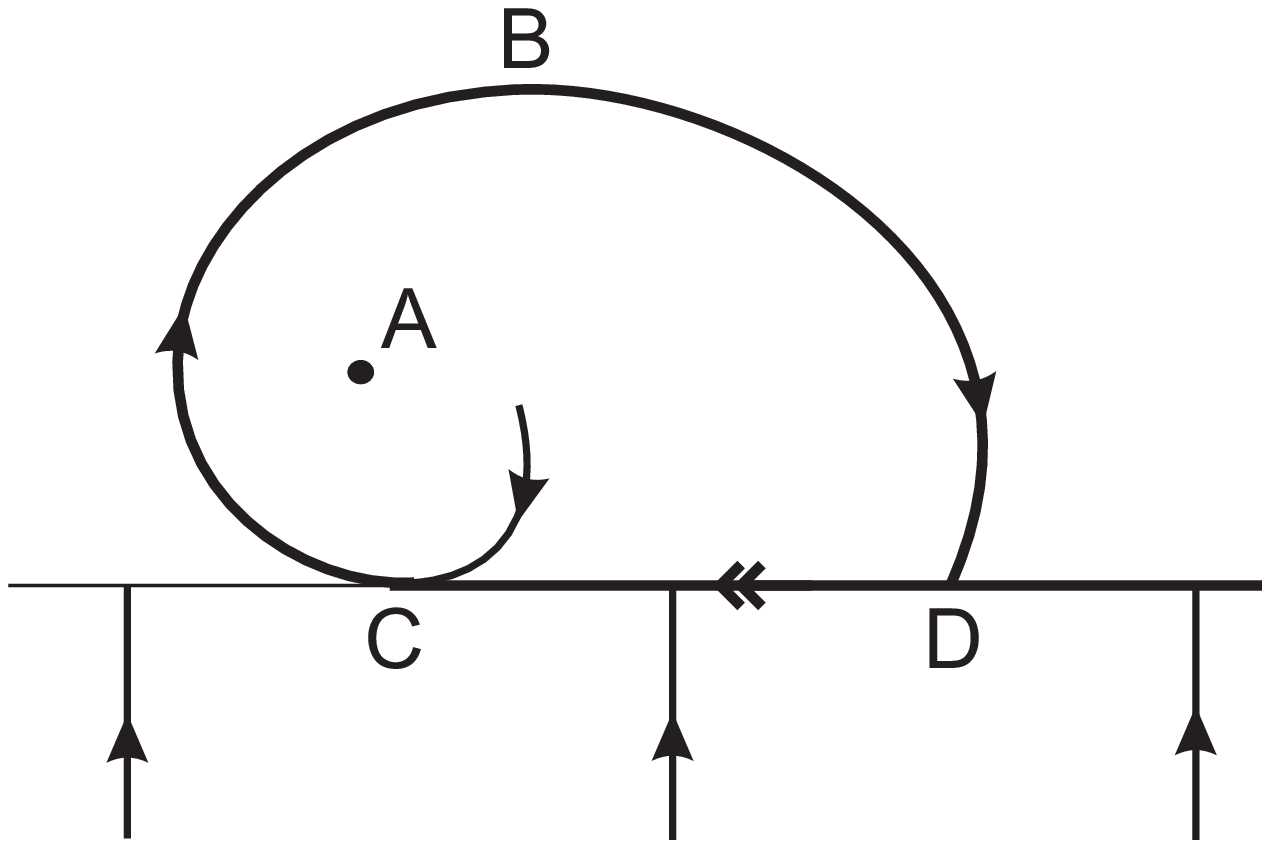}
\caption{\small{$p_1$ is repeller.}} \label{fig politrajetoria foco
atrator}
\end{center}
\end{minipage} \hfill
\begin{minipage}[b]{0.485\linewidth}
\epsfxsize=5cm \psfrag{A}{$p_{1}$} \psfrag{B}{$\gamma_{1}$}
\psfrag{C}{$A$} \psfrag{D}{$B$}
 \epsfbox{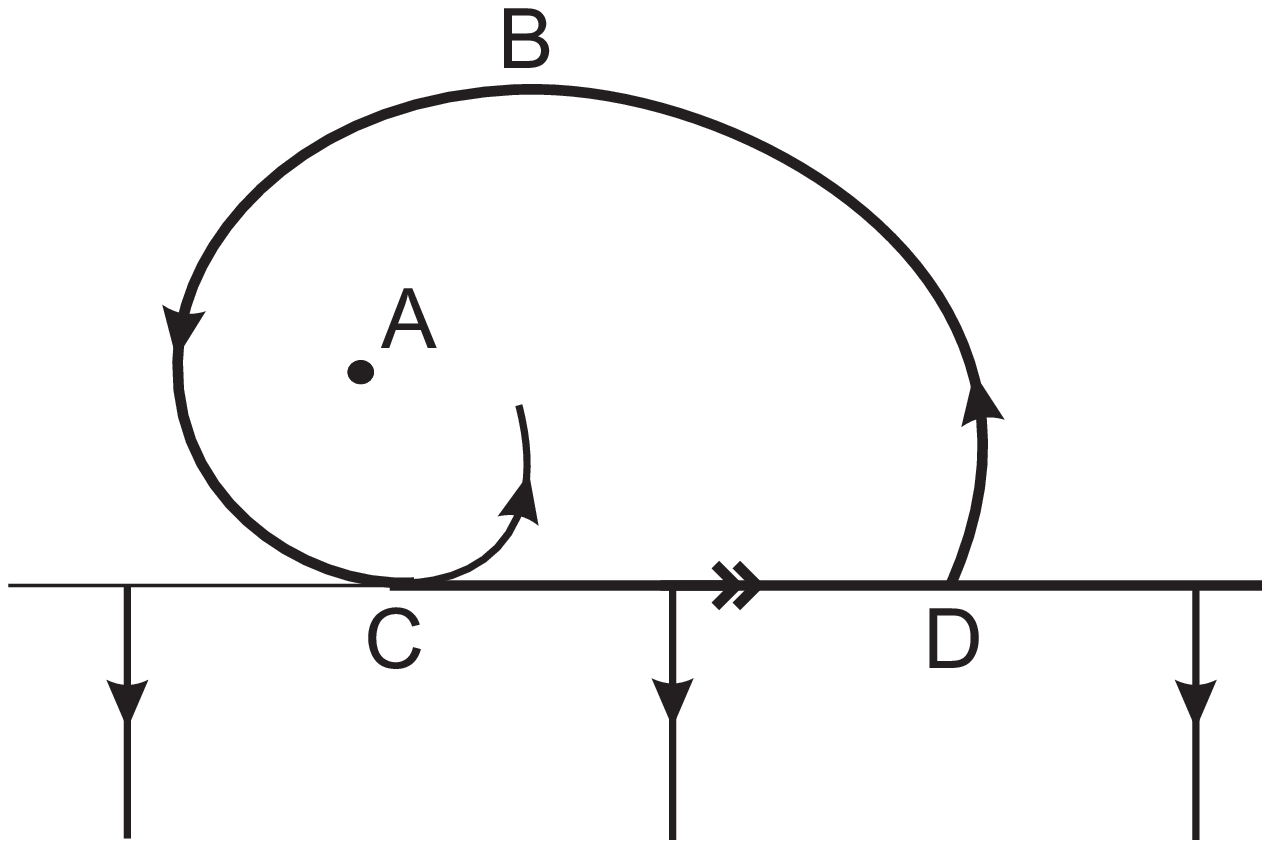}
\caption{\small{$p_1$ is attractor.}} \label{fig politrajetoria foco
repulsor}
\end{minipage}
\end{figure}

We will prove that there is not danger of ambiguity in the choose of
the closed paths, differently that what happens in the previous
example. In Figure \ref{fig politrajetoria foco atrator}, we
consider the hyperbolic canard cycle given by $\Gamma_{0} =
\gamma_{1} \cup \overline{B \, A}$. Any canard cycle of $X_0$ having
pieces of sliding region must to pass by the $\Sigma$-fold point
$A$, now walk in $\gamma_{1}$ (which is the only one possibility
that we have!) and meet the point $B$ on $\Sigma$.  The unique
choice we have is return to $A$, closing the path, without
ambiguity. For Figure \ref{fig politrajetoria foco repulsor} the
analysis is more interesting. Obviously we can use the trick of take
the vector field $-X_{0}$ and obtain an analogous result that the
previous, however we prefer give here a complete idea to the case
described in Figure \ref{fig politrajetoria foco repulsor}. Since in
the semi straight line $r=\overrightarrow{A \, B}$ we have an
escaping region, the canard cycle must have only escaping region in
its composition. Note that, if we choose to depart from $r$ by a
point in the segment $(A,B)$ then this path does not return to
$(A,B) \subset \Sigma$ (it will move spirally around the focus
$p_{1}$), if we choose going out from $r$ by a point after $B$ in
$r$ then this path also does not return to $(A,B) \subset \Sigma$
(to this path return to $(A,B)$ it must return in a sliding region,
what is not allowed because only escaping regions compose this
hyperbolic canard cycle). So we must leave $r$ by the point $B$ and
to close the curve. Therefore in any case there is not danger of
ambiguity in the choose of the closed curve and so the Poincar\'{e}
index for non-smooth vector fields is well defined. To hyperbolic
canard cycles of kind III with another particularities is enough
repeat the ideas exposed here. To hyperbolic canard cycles of kinds
I and II clearly there is not danger of ambiguity in the choose of
the closed curves once there are not escaping or sliding region in
its composition. \bbox


\begin{corollary}
Under the hypothesis of the previous theorem and assuming that all
canard cycles of $X_{0}$ are hyperbolic, we have that:
\begin{enumerate}
\item If $\Gamma_{0}$ is a canard cycle then inside $\Gamma_0$ there
exist $(2n+1)$  critical points of $X_0$, being $n$ saddles or
$\Sigma$-saddles and $(n+1)$ focus, $\Sigma$-repeller or
$\Sigma$-attractor.
\item If all critical points of $X_{0}$ are saddle or $\Sigma$-saddle then $X_{0}$ does not have canard cycles.
\end{enumerate}
\end{corollary}
\vspace{-.2cm}\noindent\textit{Proof.} Since the index of each
saddle and each $\Sigma$-saddle point is equal to $-1$ and the index
of each other critical point of $X_0$ is equal to $1$ the result is
an immediate consequence of
Theorem \ref{teoC}.\bbox\\


\section{Applications and Examples}

\subsection{Heteroclinic Orbits}

Consider the notation of the Theorem  \ref{teoA}. We give now an
example of a curve that satisfies all the hypothesis in this theorem
except that there exists a point  $C \in (A,B)$ such that the
vectors $X_{1}(C)$ and $X_{2}(C)$ are not linearly independent;
instead of $\Gamma$ obtained in the theorem we have here a
\textit{``$\Sigma$-loop"}, that is a
\textit{$\Sigma$-saddle-attractor} with connection
between $\Sigma$-separatrices.\\

\noindent\textbf{Example}. Consider the non-smooth vector field
$X_{0} = (X_{1}, X_{2})$\\ with $X_{1}(x,z) = (x+z-1,-x+z-1)$,
$X_{2}(x,z) = (-x^{2} + \frac{3}{2}x - \frac{1}{2}, 1)$ and
discontinuity set given by the $x$-axis, i.e., $f(x,z)=z$. On $z=0$,
we have $ X_{1}(x,0) = (x-1,-x-1)$ and $X_{2}(x,0) = (-x^{2} +
\frac{3}{2}x - \frac{1}{2}, 1)$ and so,
\begin{equation}
(X_{1}.f)(x,0)= -x-1 \, , \quad    (X_{2}.f)(x,0) = 1.
\end{equation}

\begin{figure}[!h]
\begin{minipage}{0.6\linewidth}
\begin{center}
\psfrag{A}{} \psfrag{B}{$-1$} \psfrag{C}{$1$}
\psfrag{D}{$\frac{1}{2}$} \psfrag{E}{$X_{1}$} \psfrag{F}{$X_{2}$}
\psfrag{G}{$\Sigma$} \epsfxsize=4.5cm \epsfbox{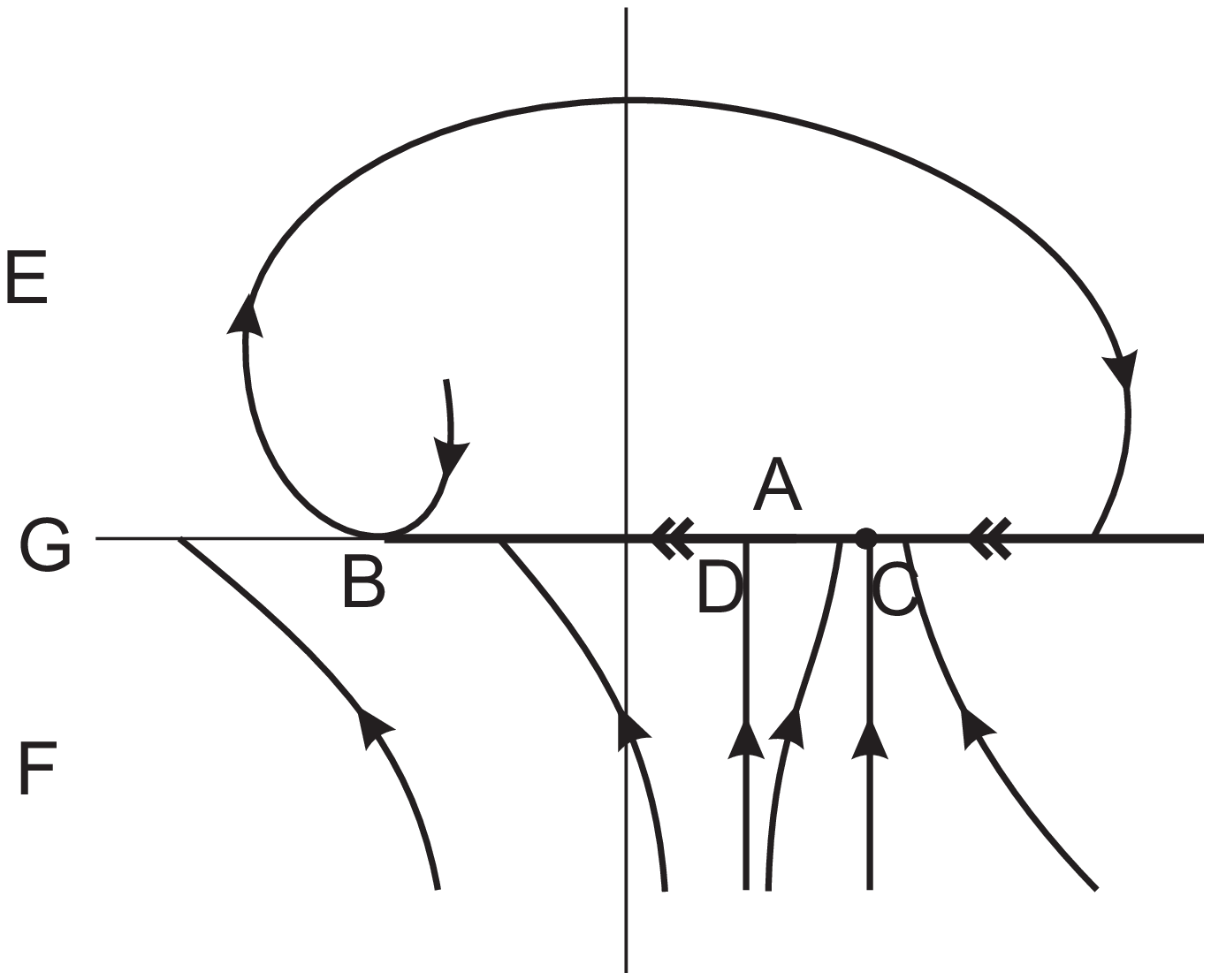}
\caption{\small{$\Sigma$-saddle-attractor with $\Sigma$-separatrices
connection.}} \label{fig d-sela-no}
\end{center}
\end{minipage} \hfill
\begin{minipage}{0.39\linewidth}
\begin{center}
\psfrag{A}{$1$} \psfrag{B}{$-1$} \psfrag{C}{$\frac{-3}{2}$}
\psfrag{H}{$h(x)$}
 \epsfxsize=4.3cm \epsfbox{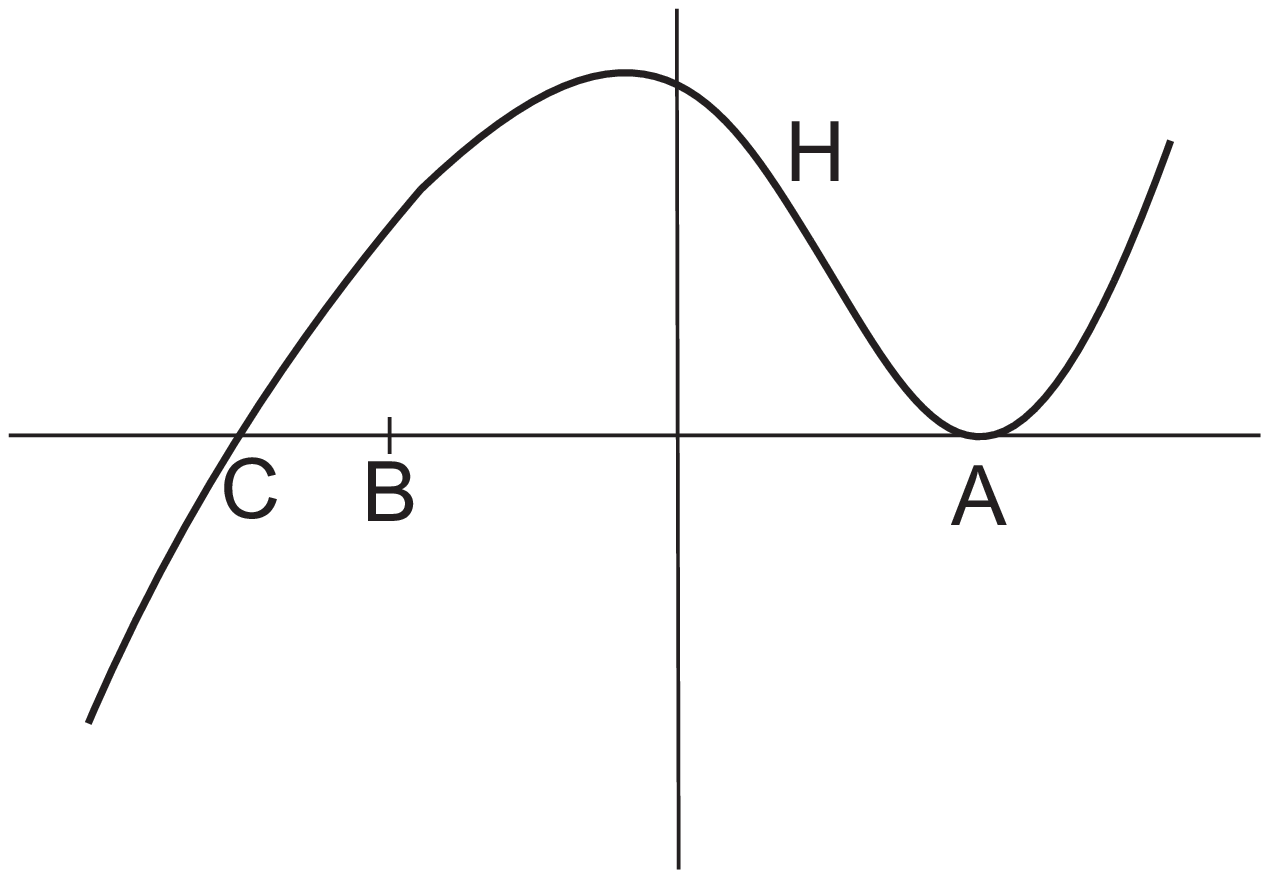}
\caption{\small{Graphic of $h$.}} \label{fig grafico h}
\end{center}
\end{minipage}
\end{figure}

In this way, we can conclude that   $x=-1$ is a $\Sigma$-fold point
of $X_{1}$ which determines a focal kind arc. For $x>-1$ we have
that $\Sigma$ is a sliding region and for $x < -1$ it is a sewing
region (see Figure \ref{fig d-sela-no}). We show now that there
exists a point $C$ in the semi straight line $x>-1$ for which $
X_{1}(C) = \lambda. X_{2}(C)$. In fact, if $X_{1}(x,0) = \lambda.
X_{2}(x,0)$ then $h(x)=(x-1)^{2}(x+ \frac{3}{2})=0$. The graphic of
$h(x)$ is given in Figure \ref{fig grafico h}. We observe that $h$
is equal to $-H$, where $H$ is the direction function defined
previously. So, we have the situation described in the Theorem
\ref{teoA}, except that $X_{1}$ and $X_{2}$ are not linearly
independent in $x=1$ where $X_0^{\Sigma}$ has an equilibrium point.
The orientation of $X_0^\Sigma$ is in direction to the $\Sigma$-fold
point because for $x=\frac{1}{2}$ we have $X_{1}(\frac{1}{2},0)=
(\frac{-1}{2},\frac{-3}{2})$, $X_{2}(\frac{1}{2},0)= (0,1)$ and so,
the direction function $H$ is negative ($H(\frac{1}{2}) = -
\frac{1}{5}$), analogously for $x=\frac{3}{2}$ we have
$X_{1}(\frac{3}{2},0)= (\frac{1}{2},\frac{-5}{2})$,
$X_{2}(\frac{1}{2},0)= (\frac{-1}{2},1)$ and so, the direction
function $H$ also is negative ($H(\frac{3}{2}) = - \frac{3}{14}$).
The  pseudo equilibrium $p=(1,0)$ is a $\Sigma$-saddle-attractor
where the $\Sigma$-separatrices are connected.

\begin{figure}[!h]
\begin{minipage}{0.52\linewidth}
\begin{center}
\psfrag{A}{$D$-singularidade} \psfrag{B}{$-1$} \psfrag{C}{$1$}
\psfrag{D}{$\frac{1}{2}$} \psfrag{E}{$X_{1}$} \psfrag{F}{$X_{2}$}
\psfrag{G}{$\Sigma$} \epsfxsize=4.5cm
\epsfbox{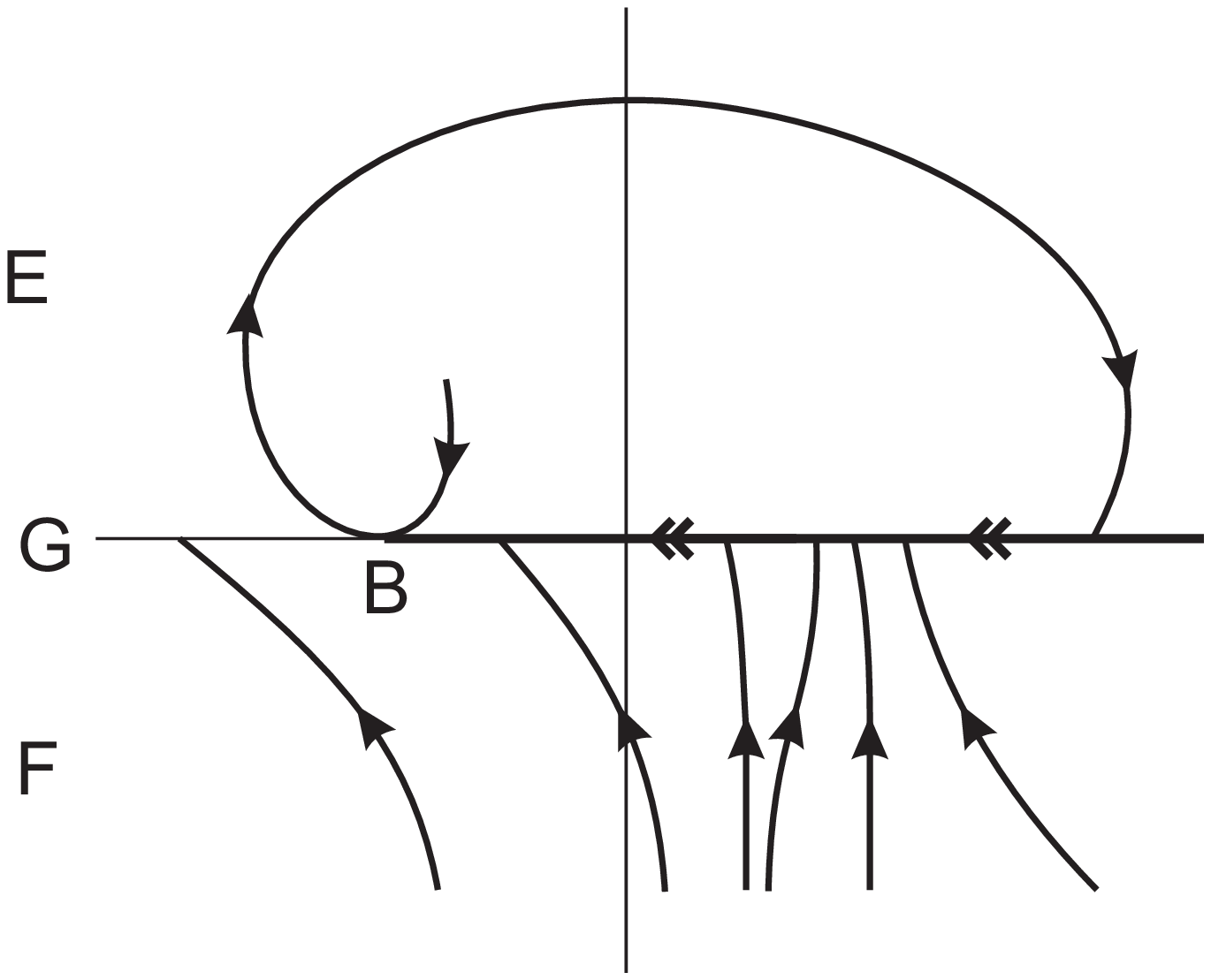} \caption{\small{Hyperbolic Canard
Cycle of kind III.}} \label{fig d-sela-no bifuracacao sem sing}
\end{center}
\end{minipage} \hfill
\begin{minipage}{0.47\linewidth}
\begin{center} \psfrag{A}{$1$} \psfrag{B}{} \psfrag{C}{$\frac{-3}{2}$}
\psfrag{H}{$h(x)$}
 \epsfxsize=4.3cm \epsfbox{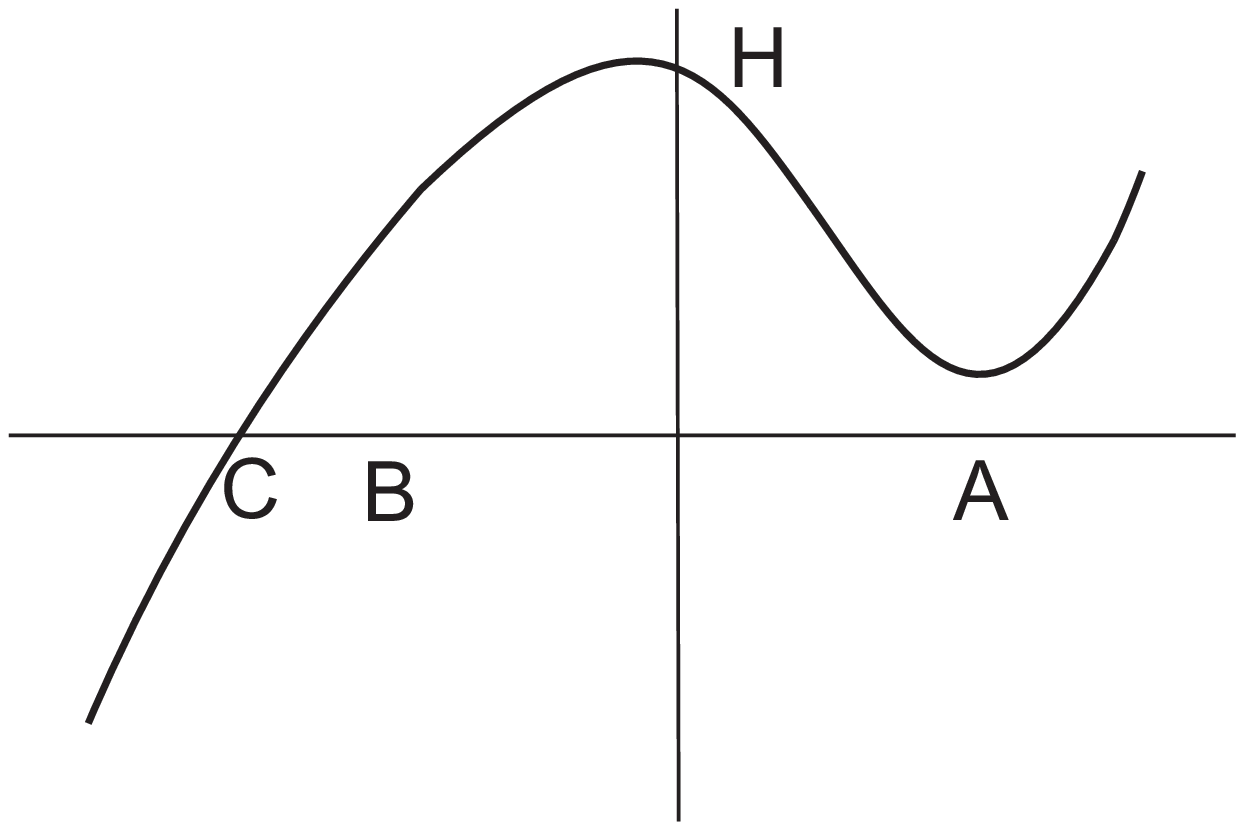}
\caption{\small{Graphic of $h$.}} \label{fig grafico h sem raiz}
\end{center}
\end{minipage}
\end{figure}

\begin{figure}[!h]
\begin{minipage}{0.52\linewidth}
\begin{center}
\psfrag{A}{$\Sigma$-saddle} \psfrag{B}{$-1$}
\psfrag{C}{$\Sigma$-attractor} \psfrag{D}{$\frac{1}{2}$}
\psfrag{E}{$X_{1}$} \psfrag{F}{$X_{2}$} \psfrag{G}{$\Sigma$}
\epsfxsize=4.5cm \epsfbox{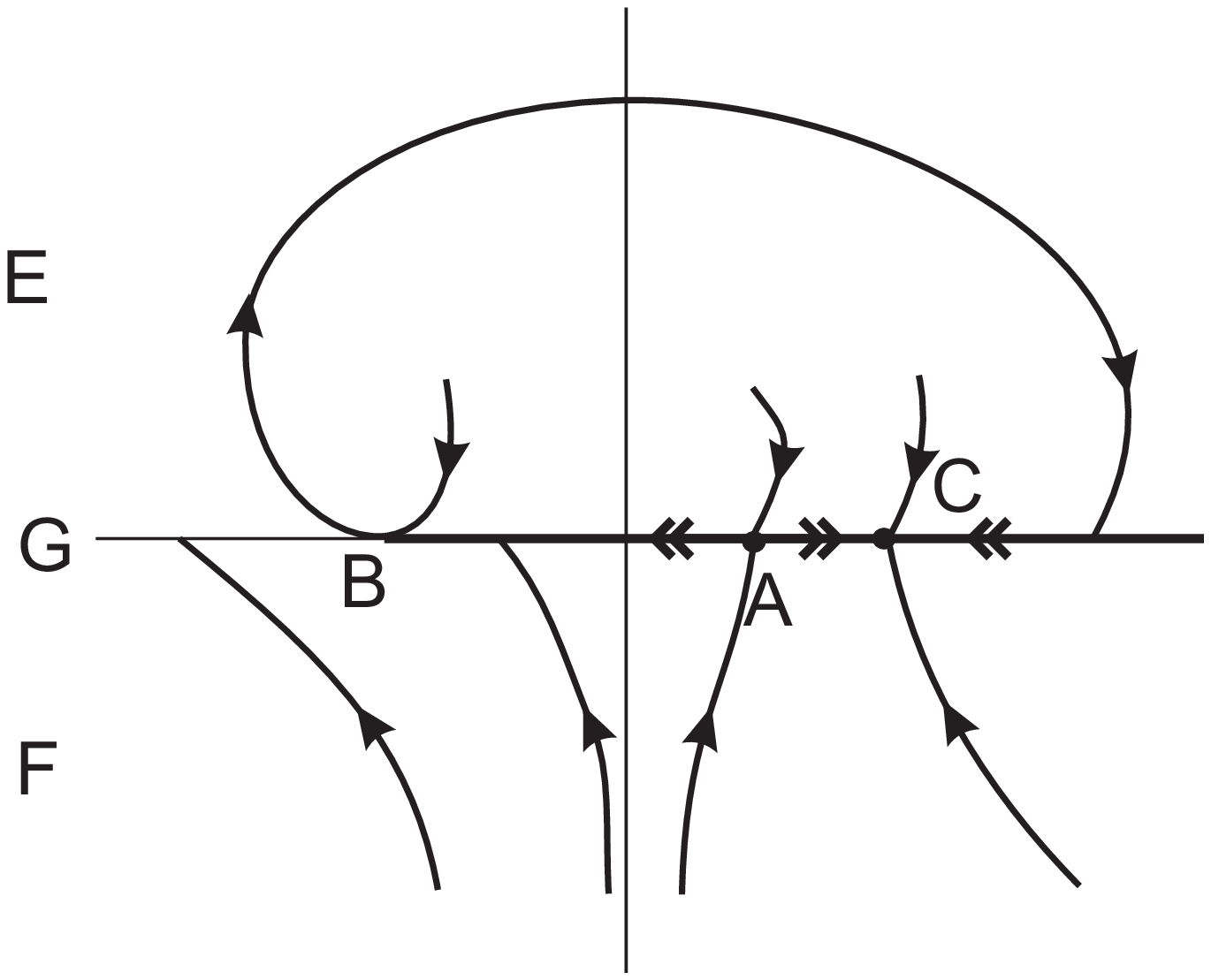}
\caption{\small{$\Sigma$-saddle and $\Sigma$-attractor with
$\Sigma$-separatrices connection.}} \label{fig d-sela-no bifurcacao
com 2 sing}
\end{center}
\end{minipage} \hfill
\begin{minipage}{0.47\linewidth}
\begin{center}\psfrag{A}{$1$} \psfrag{B}{} \psfrag{C}{$\frac{-3}{2}$}
\psfrag{H}{$h(x)$}
 \epsfxsize=4.3cm \epsfbox{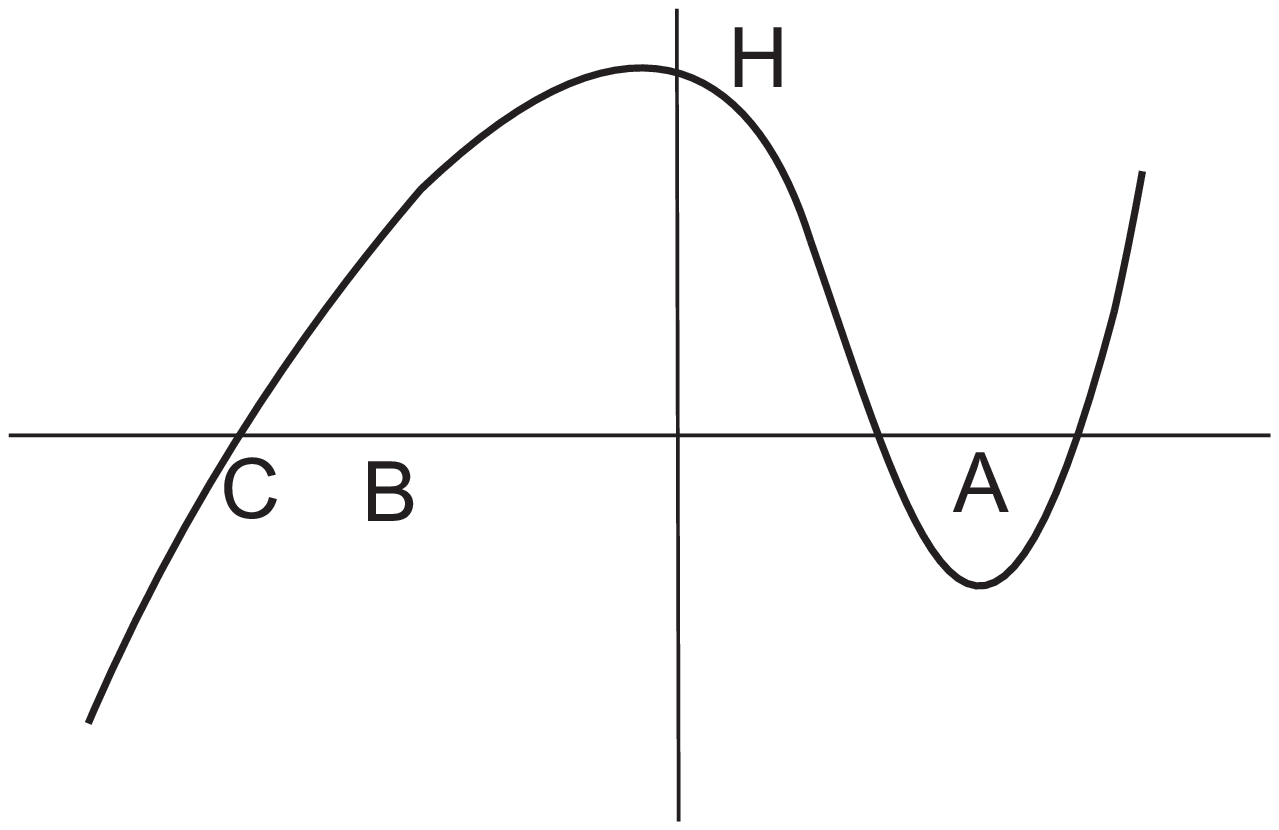}
\caption{\small{Graphic of $h$.}} \label{fig grafico h com 2 raizes}
\end{center}\end{minipage}
\end{figure}

\vspace{.5cm} \noindent\textbf{Example} \textit{(Bifurcation of the
previous example)} In the previous example, the pseudo equilibrium
with a ``loop" was found because the function $h$ has a double-zero
at $x=1$. So, we can conclude that putting small variations on the
fields $X_{1}$ and $X_{2}$ a new function $h$ appears, with two
simple real zeros (see Figure \ref{fig grafico h com 2 raizes}) or
without real zeros (see Figure \ref{fig grafico h sem raiz}) in a
neighborhood of $x=1$. The phase portrait of $X_{0}$ for small
variations on $X_{1}$ and $X_{2}$ are showed in Figures \ref{fig
d-sela-no bifuracacao sem sing} and \ref{fig d-sela-no bifurcacao
com 2 sing}. Note that for the case showed in Figure \ref{fig
d-sela-no bifuracacao sem sing}, the direction function $H$ is
always negative in a neighborhood of $x=1$ and we can apply the
corollary \ref{corteoA} to conclude that there exists a hyperbolic
canard cycle of kind III and, for the case showed in Figure \ref{fig
d-sela-no bifurcacao com 2 sing}, the direction function  $H$ has
two simple real zeros in a neighborhood of $x=1$ and assumes
positive values between this two points and negative values in the
rest of the sliding region. In this way we have a bifurcation model
where imposing small variations we can have a hyperbolic canard
cycle or we can have two pseudo equilibrium points with a stable
connection between its separatrices. We call this a
\textbf{$\Sigma$-Loop Bifurcation}.

\begin{figure}[!h]
\begin{minipage}{0.485\linewidth}
\begin{center}
\psfrag{A}{$A$} \psfrag{B}{$B$} \psfrag{C}{$\Sigma$-attractor}
\psfrag{D}{$\frac{1}{2}$} \psfrag{E}{$X_{1}$} \psfrag{F}{$X_{2}$}
\psfrag{G}{$\Sigma$} \epsfxsize=4.5cm
\epsfbox{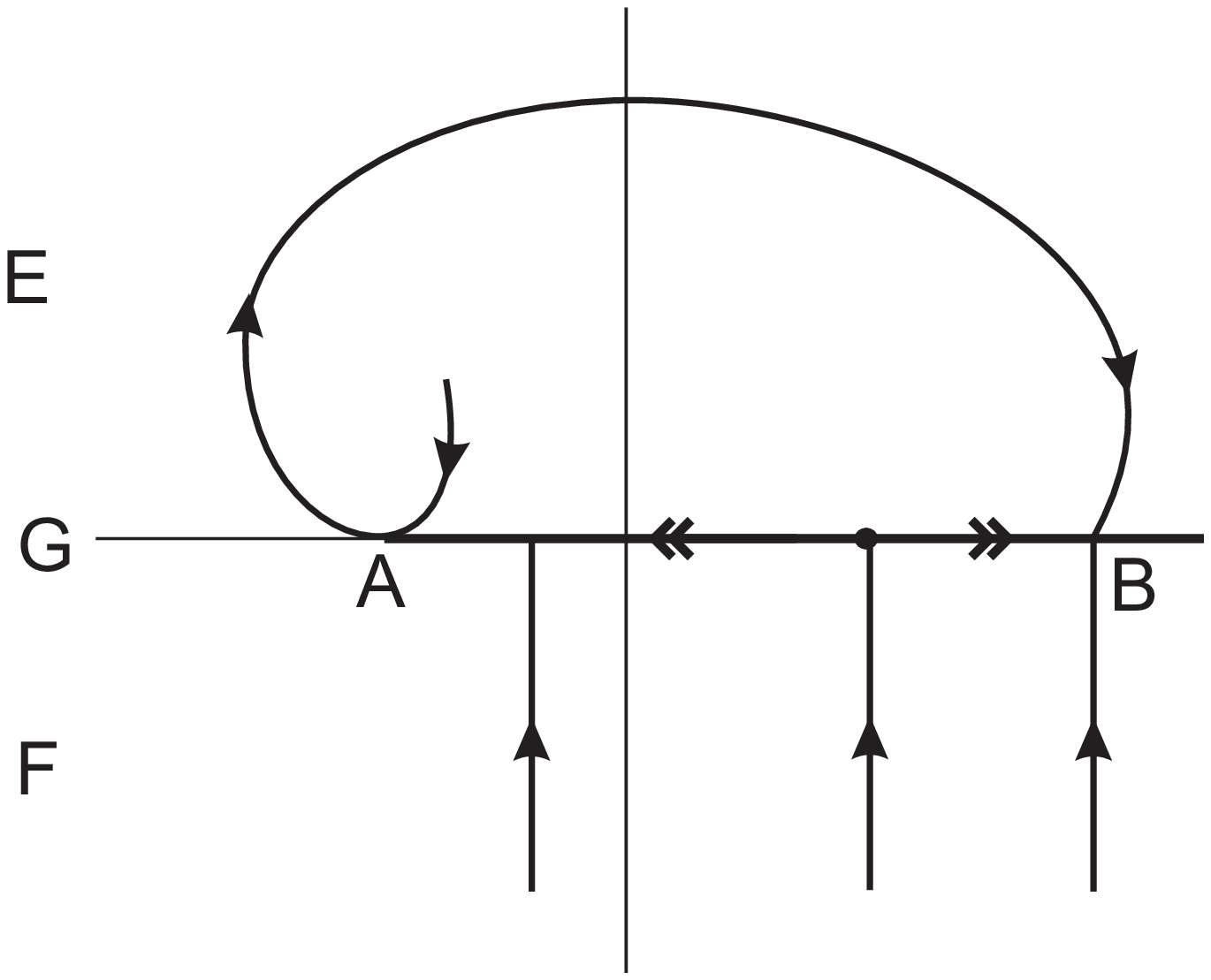} \caption{\small{Stable
configuration.}} \label{fig d-loop repulsor}
\end{center}
\end{minipage} \hfill
\begin{minipage}{0.485\linewidth}
\begin{center}\psfrag{A}{$A$} \psfrag{B}{$B$} \psfrag{C}{$\Sigma$-attractor}
\psfrag{D}{$\frac{1}{2}$} \psfrag{E}{$X_{1}$} \psfrag{F}{$X_{2}$}
\psfrag{G}{$\Sigma$}
 \epsfxsize=4.5cm \epsfbox{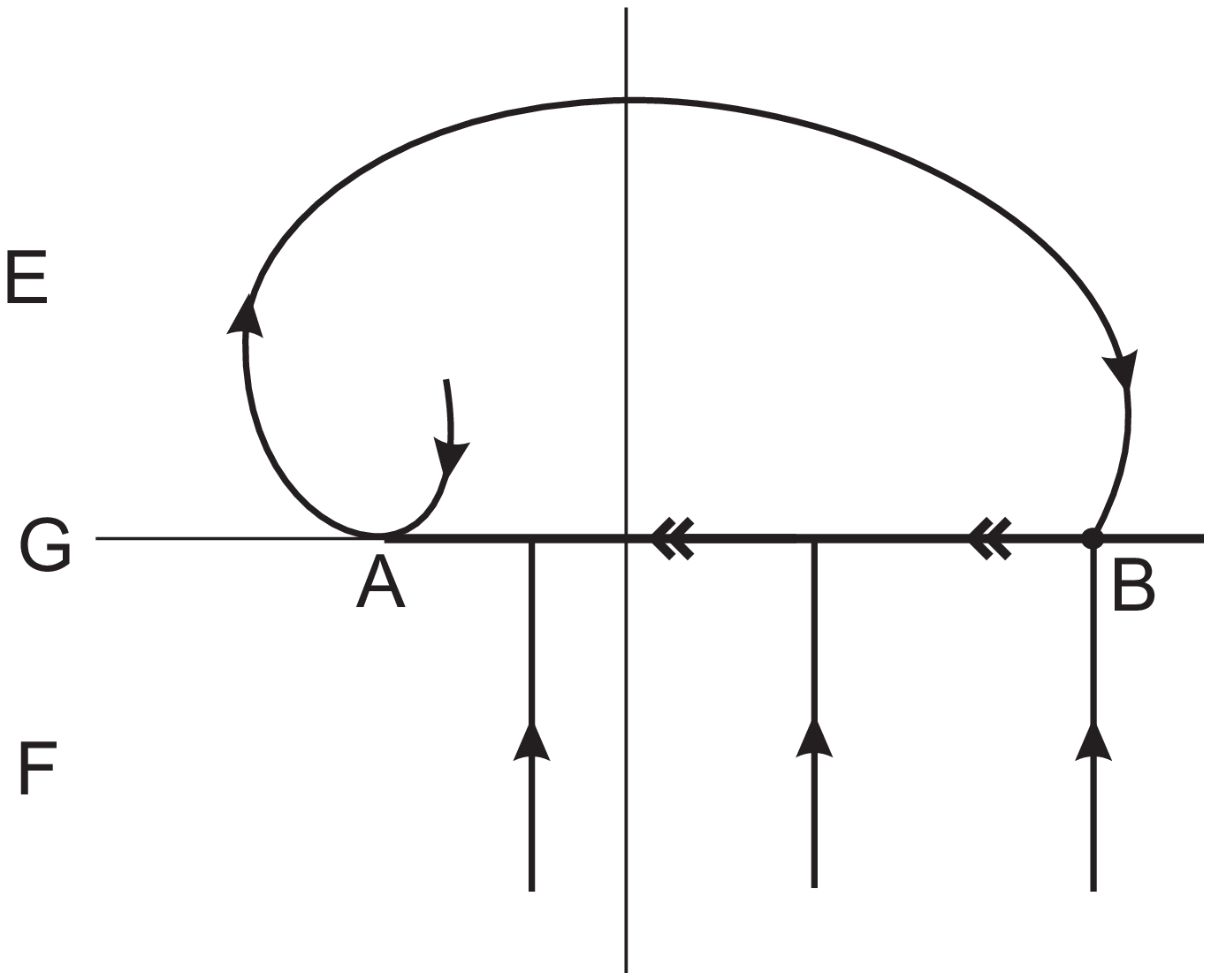}
\caption{\small{Unstable configuration.}} \label{fig d-loop
diferente}
\end{center}\end{minipage}
\end{figure}

Following the notation of Theorem \ref{teoA} and the ideas just
exposed we can state the next Proposition:

\begin{proposition}Let $X_0$ a non-smooth vector field.
\begin{enumerate}
\item $X_0$ has an unstable configuration topologically equivalent to that one in Figure \ref{fig d-sela-no}
if and only if (i) the direction function $H$ is well defined in
$[A,B]$, (ii) $H$ has a single zero in $(A,B)$ and (iii) $H(B) < 0$.

\item $X_0$ has a stable configuration topologically
equivalent to that one in Figure \ref{fig d-loop repulsor} if and
only if (i) the direction function $H$ is well defined in $[A,B]$,
(ii) $H$ has a single zero in $(A,B)$ and (iii) $H(B) > 0$.

\item $X_0$ has
an unstable configuration topologically equivalent to that one in
Figure \ref{fig d-loop diferente} if and only if (i) the direction
function $H$ is well defined in $[A,B]$, (ii) $H$ do not have zeros
in $(A,B)$ and (iii) $H(B)= 0$.

\end{enumerate}
\end{proposition}

\vspace{-.25cm} \noindent\textit{Proof.} It is straightforward
following what is done in the previous example.\\

Moreover, concerning with the item (1) of the Proposition 1 small
perturbations in $X_0$ produces the effects showed in the previous
example and  concerning with the item (3) of the Proposition 1 small
perturbations in $X_0$ produces  effects described in Theorem
\ref{teoA}, in the previous example and in the item (2) of the
Proposition 1.

\subsection{Canard Cycles and Singular Perturbations Problems.} In this section we show how
the regularization process gives a singular perturbation problem. In
this context, the canard cycles defined in this paper can be
considered as \emph{limit periodic
sets} of singular problems. First of all we present some basic definitions.\\

\vspace{-.3cm}

\begin{definition} Let $U\subseteq \R^2$ be an open subset and take $\e\geqslant 0$.
A singular perturbation problem in $U$ (SP--Problem) is a
differential system which can be  written like
\begin{equation}
\label{fast} x'=dx/d\tau=l(x,y,\e),\quad  y'=dy/d\tau=\e m(x,y,\e)
\end{equation} or equivalently, after the time re-scaling $t=\e\tau$
\begin{equation}
\label{slow} \e{\dot x}=\e dx/dt=l(x,y,\e),\quad {\dot y}=dy/dt=
m(x,y,\e),
\end{equation}
with $(x,y)\in U$ and $l,m$ smooth in all variables.
\end{definition}

The understanding of the  phase portrait of the vector field
associated to a SP-problem is the main goal of the \textit{geometric
singular perturbation theory} (GSP-theory). The techniques of
GSP-theory can be used to obtain information on the dynamics of
(\ref{fast}) for small values of $\e>0,$ mainly in searching limit
cycles. System (\ref{fast}) is called the \textit{fast system}, and
(\ref{slow}) the \textit{slow system} of SP-problem. Observe that
for $\e >0$ the phase portraits of the fast and the slow systems
coincide. For $\epsilon =0,$ let $\mathcal{S}$ be the set
\begin{equation}
 \label{SM} \mathcal{S}=\left\lbrace (x,y):f(x,y,0)=0\right\rbrace
 \end{equation}
of all singular points of (\ref{fast}). We call $\mathcal{S}$ the
slow manifold of the singular perturbation problem and it is
important to notice that equation (\ref{slow}) defines  a dynamical
system, on $\mathcal{S}$, called the {\it reduced problem}:
\begin{equation}
\label{reduced} f(x,y,0)=0 ,\quad {\dot y}= g(x,y,0).
\end{equation}
Combining results on the dynamics of these two limiting problems,
with $\e = 0$, one obtains information on the dynamics of $X_\e$ for
small values of $\e$. We refer to \cite{F} for an introduction to
the general theory of singular perturbations. Related problems can
be seen in \cite{BST}, \cite{DR} and \cite{S}. Let us apply the
techniques of GSP-Theory to study  hyperbolic canard cycles. \\

\noindent\textbf{Example.} Consider the non-smooth vector field $X_0
= (X_1,X_2)$ with $X_{1}(x,y)=(x+y-1,-x+y+1)$, $X_{2}(x,y) = (1,2)$
and $f(x,y) = x$. The regularized vector field becomes
$$\begin{array}{lcl}
   \dot{x} & = & \frac{x+y}{2} +
\varphi \left( \frac{x}{\epsilon} \right)\frac{x+y-2}{2} \, ,
\\
   \dot{y}  & =  & \frac{-x+y+3}{2} +
\varphi \left( \frac{x}{\epsilon}
\right)\frac{-x+y-1}{2}.\end{array}
$$
where $\varphi(\frac{x}{\epsilon})$ is the transition function.
Making the change of variables  $x=r. \cos \theta$ and $\epsilon =
r. \sin \theta$ we obtain
\begin{equation}\label{eq sp problem com mudanca de variaveis}
\begin{array}{lcl}
  r \dot{\theta}  & = & -\sin \theta \left(  \frac{r.\cos\theta + y}{2} + \varphi(\cot \theta)\frac{r.\cos\theta + y -2}{2} \right) \, , \\
  \dot{y} &= & \frac{-r.\cos\theta+y + 3}{2} + \varphi(\cot
\theta)\frac{-r.\cos\theta +y - 1}{2}.
\end{array}
\end{equation}

In the blowing up locus $r=0$ the fast dynamics is determined by the
system
 \[ \theta'   =  -\sin \theta \left( \frac{y}{2} + \varphi(\cot \theta)\frac{y -2}{2} \right) \, , \quad y'=0
\, ;  \] and the slow dynamics on the slow manifold is determined by
the reduced system
\[  \frac{y}{2} + \varphi(\cot \theta)\frac{y -2}{2} = 0 \, , \quad  \dot{y}= \frac{y + 3}{2} + \varphi(\cot
\theta)\frac{y-1}{2}.
\]

We remark that the slow manifold is implicitly defined by
$\frac{y}{2} + \varphi(\cot \theta)\frac{y -2}{2} = 0$ and
$y(\theta)$ defined in this way is such that $
\displaystyle\lim_{\theta \longrightarrow \frac{\pi}{4}}y(\theta) =
1 $, $ \displaystyle\lim_{\theta \longrightarrow
\frac{3\pi}{4}}y(\theta) = -\infty $.

\begin{figure}[!h]
\begin{minipage}{0.485\linewidth}
\psfrag{A}{$X_{2}$} \psfrag{B}{$X_{1}$} \psfrag{C}{$y$}
\psfrag{D}{$\pi$} \psfrag{E}{$\frac{3 \pi}{4}$}
\psfrag{F}{$\frac{\pi}{2}$}
\psfrag{G}{$\frac{\pi}{4}$}\psfrag{H}{$0$}\psfrag{I}{$x$}
\epsfxsize=6cm \epsfbox{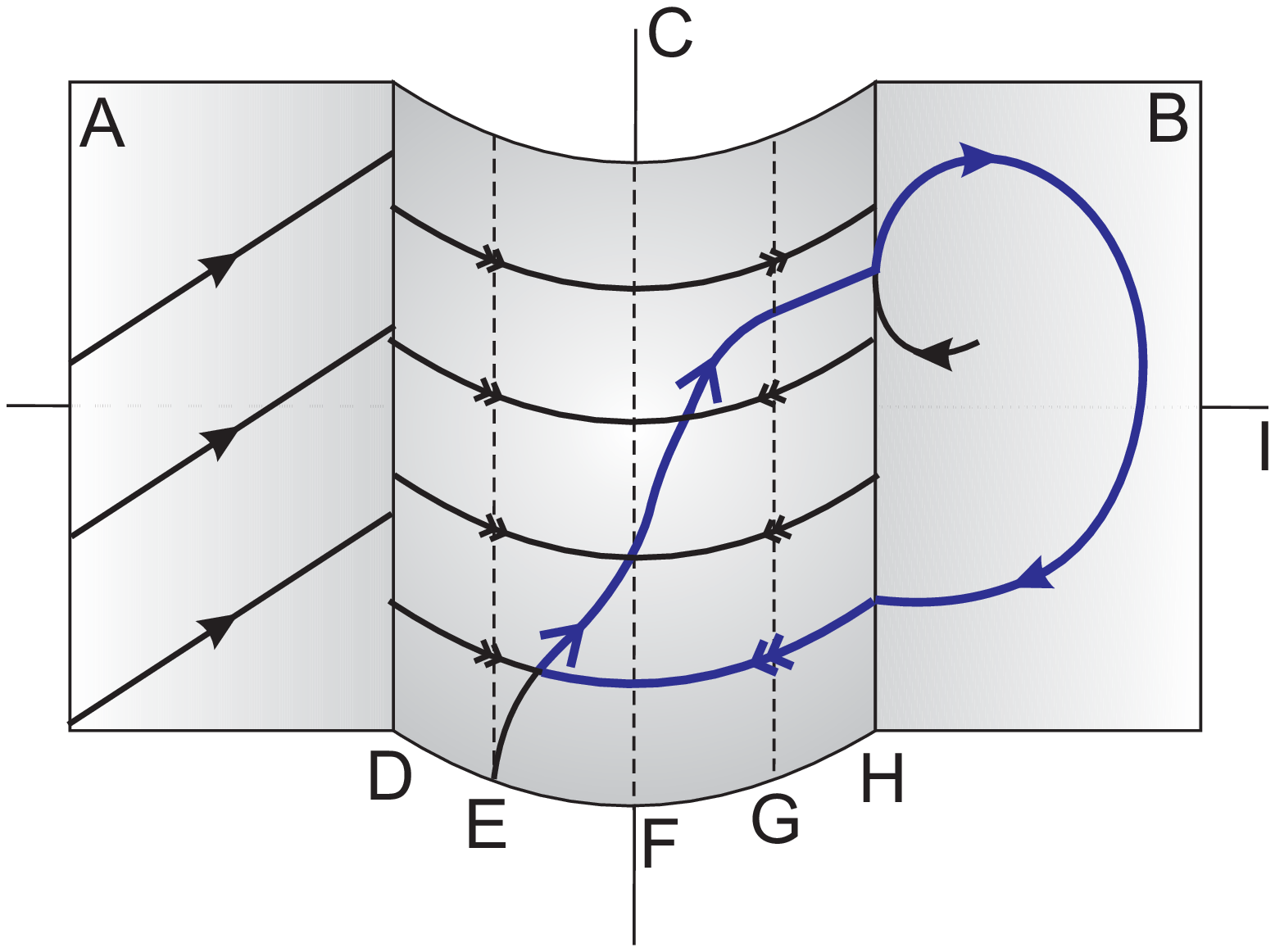}
\caption{\small{Singular perturbation of a Canard with one
$\Sigma$-fold point.}} \label{fig pert singular final}
\end{minipage} \hfill
\begin{minipage}{0.485\linewidth}
\psfrag{C}{$y$}  \psfrag{D}{$x$} \psfrag{A}{$X_{2}$}
\psfrag{B}{$X_{1}$} \psfrag{F}{$\frac{\pi}{2}$} \epsfxsize=6.2cm
\epsfbox{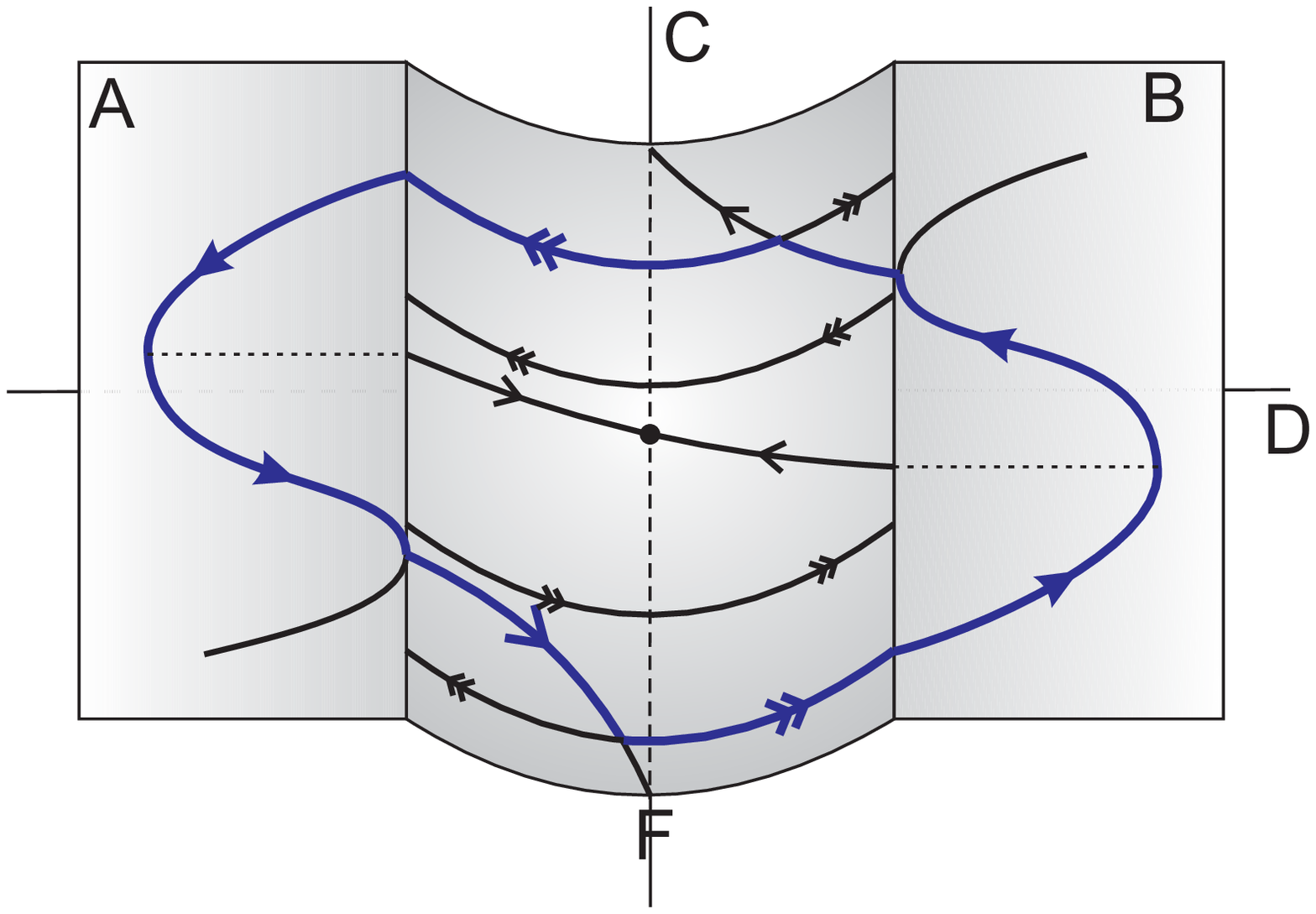} \caption{\small{Singular
perturbation of a Canard with two $\Sigma$-fold points.}} \label{fig
pert singular exemplo deles}
\end{minipage}
\end{figure}

In the phase portrait on the blowing up locus double arrow over one
the trajectory means that the trajectory is of the fast dynamical
system, and simple arrow means that the trajectory is of the slow
dynamical system. So, we can draw the slow variety and its
orientation and give the orientation of the fast flow (see Figure
\ref{fig pert singular final}).


In the final example of \cite{Claudio-PR-Marco} the authors apply
the GSP-Theory to the non-smooth vector field $X_0(x,y) =
(X_1(x,y),X_2(x,y)) = ((3y^2-y-2,1),(-3y^2-y+2,-1))$ and obtain the
SP-problem which behavior is described in Figure \ref{fig pert
singular exemplo deles}.\\

%
%
%
%
%

In \cite{LST} the authors prove that $\Sigma^2 \cap \Sigma^3$ is
homeomorphic to the slow variety and that the sliding vector field
$X^{\Sigma}_{0}$ is topologically equivalent to the reduced problem.
So, we can apply step-by-step the method described in section 4 of
this paper and found  tubular neighborhoods for the canard cycles in
Figures \ref{fig pert singular final} and \ref{fig pert singular
exemplo deles}. Moreover, we can apply the Theorem \ref{teoB} and
conclude that it is a limit set (making $r \rightarrow 0$) of
hyperbolic limit cycles. This also is true to any one hyperbolic
canard cycle.\\

\noindent {\textbf{Acknowledgments.}} The first and the third
authors are partially supported by CAPES and CNPq, the second author
is partially supported by a FAPESP-BRAZIL grant 2007/08707-5.


\begin{thebibliography}{99}

%

\bibitem{BST} {\sc C.A. Buzzi, P.R. da Silva and M.A. Teixeira},
{\it Singular Perturbation Problems For Time Reversible Systems},
Proc. Amer. Math. Soc., \textbf{133} (2005), 3323-3331.

\bibitem{Claudio-PR-Marco} {\sc C.A. Buzzi, P.R. da Silva  and M. A. Teixeira},
{\it A Singular Approach To Discontinuous Vector Fields On The
Plane}, Journal of Differential Equations, \textbf{231} (2006),
633-655.
%

\bibitem{DR} {\sc F. Dumortier and R. Roussarie},
{\it Canard Cycles And Center Manifolds}, Memoirs Amer. Mat. Soc.
\textbf{121}, 1996.

\bibitem{F} {\sc N. Fenichel},
{\it Geometric Singular Perturbation Theory For Ordinary
Differential Equations}, Journal of Differential Equations
\textbf{31} (1979), 53--98.

\bibitem{Fi} {\sc A.F. Filippov},
{\it Differential Equations With Discontinuous Righthand Sides},
Mathematics and its Applications (Soviet Series), Kluwer Academic
Publishers-Dordrecht, 1988.




\bibitem{K} {\sc V. S. Kozlova},
{\it Roughness Of A Discontinuous System}, Vestinik Moskovskogo
Universiteta, Matematika \textbf{5} (1984), 16--20.



\bibitem{LST} {\sc J. Llibre, P.R. Silva and M.A. Teixeira},
{\it  Sliding Vector Fields Via Slow-Fast Systems}, Bulletin of the
Belgian Mathematical Society Simon Stevin \textbf{15-5} (2008),
851--869.



\bibitem{SM} {\sc J. Sotomayor and A.L. Machado},
{\it  Structurally Stable Discontinuous Vector Fields On The Plane},
Qual. Theory of Dynamical Systems, \textbf{3} (2002), 227--250.



\bibitem{ST} {\sc J. Sotomayor and M.A. Teixeira},
{\it Regularization Of  Discontinuous Vector Fields}, International
Conference on Differential Equations, Lisboa (1996), 207--223.

\bibitem{SZ} {\sc J. Sotomayor and M. Zhitomirskii}, {\it Impasse
Singularities of Differential Systems of the Form $A(x)x'=F(x)$}, J.
Diff. Equations \textbf{169}, n$^{o}$2, (2001), 567--587. 


\bibitem{S} {\sc P. Szmolyan},
{\it Transversal Heteroclinic And Homoclinic Orbits In Singular
Perturbation Problems}, Journal of Differential Equations
\textbf{92} (1991), 252--281.




\bibitem{T} {\sc M.A. Teixeira},
{\it  Generic Singularities Of Discontinuous Vector Fields}, An. Ac.
Bras. Cienc. \textbf{53}, n$^{o}$2, (1991), 257--260.

\end{thebibliography}
\end{document}